\documentclass{conm-p-l}

\usepackage{tikz-cd,url}
\usepackage{hyperref}

\usepackage{amssymb}

\usepackage{graphicx}

\usepackage[cmtip,all]{xy}

\usepackage{}

\newtheorem{theorem}{Theorem}

\newtheorem{proposition}{Proposition}

\theoremstyle{definition}

\theoremstyle{remark}
\newtheorem{remark}{Remark}

\begin{document}

 \title[Riemann-Roch spaces in the Hilbert class field]{Explicit Riemann-Roch spaces in the Hilbert class field}


\author{Jean-Marc Couveignes}
\address{Univ. Bordeaux, CNRS, INRIA, 
  Bordeaux-INP, IMB, UMR 5251, F-33400 Talence, France.}
\email{jean-marc.couveignes@math.u-bordeaux.fr}

\author{Jean Gasnier}
\address{Univ. Bordeaux, CNRS, INRIA, 
  Bordeaux-INP, IMB, UMR 5251, F-33400 Talence, France.}
\email{jean.gasnier@math.u-bordeaux.fr}


\date{May 22, 2024}

\newcommand\ddaleth{{\varrho}}
\newcommand\taum{{\tau_{{\mathrm{max}}}}}
\newcommand\dbasic{{d_{{\mathrm{basic}}}}}
\newcommand\deltabasic{{\delta_{{\mathrm{basic}}}}}
\newcommand\Qm{{Q_{{\mathrm{max}}}}}
\newcommand\Ym{{Y_{{\mathrm{max}}}}}
\newcommand\bK{{\mathbf K}}
\newcommand\bff{{\mathbf{fr}}}
\newcommand\bF{{\mathbf F}}
\newcommand\bC{{\mathbf C}}
\newcommand\itOmega{{\mathit \Omega}}
\newcommand\bA{{\mathbf A}}
\newcommand\bR{{\mathbf R}}
\newcommand\bk{{\mathbf k}}
\newcommand\bZ{{\mathbf Z}}
\newcommand\bN{{\mathbf N}}
\newcommand\bKb{{\bar {\mathbf K}}}
\newcommand\bL{{\mathbf L}}
\newcommand\up{\mathop{{\uparrow}}\nolimits }
\newcommand\Supp{\mathop{{\mathrm{Supp}}}\nolimits }
\newcommand\Epsilon{{\Supp (\epsilon)}}
\newcommand\down{\mathop{{\downarrow}}\nolimits }
\newcommand\Jac{\mathop{\rm Jac }}
\newcommand\cS{\mathcal S}
\newcommand\co{\mathfrak o}
\newcommand\cD{\mathcal D}
\newcommand\cP{\mathcal P}
\newcommand\cM{\mathcal M}
\newcommand\cR{\mathcal R}
\newcommand\cA{\mathcal A}
\newcommand\cE{\mathcal E}
\newcommand\cC{\mathcal C}
\newcommand\cW{\mathcal W}
\newcommand\cZ{\mathcal Z}
\newcommand\cU{\mathcal U}
\newcommand\cI{\mathcal I}
\newcommand\cL{\mathcal L}
\newcommand\cO{\mathcal O}
\newcommand\cB{\mathcal B}
\newcommand\Aut{\mathop{\rm{Aut}}\nolimits }
\newcommand\Res{\mathop{\rm{res}}\nolimits }
\newcommand{\mmu}{\mbox{$\raisebox{-0.59ex}
  {$l$}\hspace{-0.18em}\mu\hspace{-0.88em}\raisebox{-0.98ex}{\scalebox{2}
  {$\color{white}.$}}\hspace{-0.416em}\raisebox{+0.88ex}
  {$\color{white}.$}\hspace{0.46em}$}{}}
\newcommand\rank{\mathop{\rm{rank}}\nolimits }
\newcommand\Proj{\mathop{\rm{Proj}}\nolimits }
\newcommand\Hom{\mathop{\rm{Hom}}\nolimits }
\newcommand\HomS{\mathop{\rm{Hom}_{\rm\bf set}}\nolimits }
\newcommand\Bil{\mathop{\rm{Bil}}\nolimits }
\newcommand\Ker{\mathop{\rm{Ker}}\nolimits }
\newcommand\Gal{\mathop{\rm{Gal}}\nolimits }
\newcommand\PGL{\mathop{\rm{PGL}}\nolimits }
\newcommand\Pic{\mathop{\rm{Pic}}\nolimits }
\newcommand\Div{\mathop{\rm{Div}}\nolimits }
\newcommand\Spec{\mathop{\rm{Spec}}\nolimits }
\newcommand\Tr{\mathop{\rm{Tr}}\nolimits }
\newcommand\GL{\mathop{\rm{GL}}\nolimits }
\newcommand\FT{\mathop{\rm{FT}}\nolimits }
\newcommand\End{\mathop{\rm End }}
\newcommand\cQ{{\mathcal Q}}

\begin{abstract}
  Let $\bK$ be a finite field, $X$ and $Y$ two
  curves over $\bK$, and  $Y\rightarrow X$ an unramified abelian
  cover with Galois group $G$. Let $D$ be a divisor on $X$
  and  $E$ its pullback on $Y$. Under mild conditions the linear space
  associated with $E$ is a free $\bK[G]$-module. We study the
  algorithmic aspects and applications of these modules.
\end{abstract}

\maketitle

\section{Introduction}

Given   a curve $Y$ over a field $\bK$, and two divisors
$E$ and $Q$ on $Y$, with $Q$ effective and disjoint from $E$, the
evaluation  map
$e : H^0(Y,\cO_Y (E))\rightarrow H^0 (Y,\cO_Y/\cO_Y(-Q))$ is a natural
$\bK$-linear datum of some importance for various
algorithmic problems such  as efficient computing in the Picard group
of $Y$ (see \cite{mak1,mak2}), constructing good error correcting codes
\cite{gop1,gop2,tvc}, 
or bounding the bilinear complexity of multiplication
in finite fields \cite{STV,SH,BR,BAL,CHA,RAN}.
Assume that  $G$ is a finite group of automorphisms of $Y/\bK$, and
the divisors $E$ and $Q$ are $G$-equivariant (they are equal to their pullback by any element of $G$).
The evaluation map $e$ is then a $\bK[G]$-linear map between two $\bK [G]$-modules.
In some cases these modules  can be shown to  be  both  free. Their rank as $\bK [G]$-modules is then smaller than
their dimension as $\bK$-vector spaces, by a factor $\co$, the order of $G$. This is helpful  when $G$ is commutative,  because multiplication in $\bK [G]$ is achieved in quasi-linear time
using a discrete Fourier transform, and the advantage  of lowering dimension is stronger than
the disadvantage  of dealing with a larger  ring of scalars. 
We will
focus on free $\bK[G]$-modules arising from commutative groups acting freely
on a curve.
This special case
has a rich mathematical background and produces interesting constructions.
For example
Theorem \ref{th} states  the existence of excellent algebraic geometry codes
that can be encoded in quasi-linear time and decoded in quasi-quadratic time in their length. 

In Section \ref{sec:duakg} we review elementary properties
of $\bK[G]$-modules when $\bK$ is a commutative field
and $G$ a finite group. 
We recall in Section \ref{sec:curveact} how unramified
fibers of Galois covers of curves  produce
free $\bK[G]$-modules and we introduce natural
bases for these modules. 
We study the abelian unramified case in Section \ref{sec:comm}.
Theorem \ref{prop:free} states 
that  in this case, the Riemann-Roch space associated to a $G$-equivariant
divisor of large enough degree is a  free $\bK[G]$-module.
Evaluating at another $G$-equivariant divisor then produces
a $\bK[G]$-linear map between two free $\bK[G]$-modules.
This makes it possible to treat evaluation and interpolation
as $\bK[G]$-linear problems. We introduce the matrices associated to these problems.
Section \ref{sec:pade} is devoted to the definition and computation of Pad{\'e} approximants in this
context. The complexity of arithmetic
operations in $\bK [G]$ is bounded in Section \ref{sec:ga} using various classical discrete
Fourier transforms. Theorem \ref{prop:fifi} states that the complexity of multiplication in
$\bK [G]$ is quasi-linear when $G$ is commutative.
In Section \ref{sec:constf} we use effective 
class field theory and the algorithmics of
curves and jacobian varieties to compute
the evaluation and interpolation matrices 
introduced in Section \ref{sec:comm}.
Section \ref{sec:interpol} is concerned with  two applications
of interpolation with  $\bK[G]$-modules: multiplication
in finite fields and geometric codes.
The asymptotic properties of the  codes constructed in  this way
are studied in Section \ref{sec:gc}. We thank the anonymous referee for their comments and suggestions.
The calculation in Section
\ref{sec:excftj}
has been implemented
using 
SageMath (Version 9.4) the Sage Mathematics Software System \cite{sagemath}.

\vskip -1cm

\tableofcontents

\section{Duality for $\bK[G]$-modules}\label{sec:duakg}

In this section  $\bK$ is a commutative field and $G$
is a finite group. We state elementary properties of
$\bK[G]$-modules  and their duals. In Section~\ref{sec:inf}
we describe the natural correspondence between
$G$-invariant $\bK$-bilinear forms and $\bK[G]$-bilinear forms.
We  see in Section \ref{sec:orth} that the orthogonal
of a $\bK[G]$-submodule for either form  is the same.
Sections \ref{sec:duamo}  is concerned with the
canonical bilinear form relating a $\bK[G]$-module and
its dual.
The ring $\bK[G]$ has the Frobenius property
\cite[Chapter IX]{currei}. We recall 
in Section \ref{sec:freedi} a convenient consequence of
it.

\subsection{Invariant bilinear forms}\label{sec:inf}

Let $M$ be a right $\bK[G]$-module.
Let $N$ be a left  $\bK[G]$-module.
Let \[<.,.> \, : M\times N \rightarrow \bK \]
be a $\bK$-bilinear form. We assume that this form is invariant
by the action of $G$ in the sense that
\[<m.\sigma , n>\, = \, <m, \sigma . n >\] for
every $m$ in $M$,  $n$ in $N$, and $\sigma$ in $G$.
We define a map
  \begin{equation}\label{eq:inf} 
\xymatrix@R-2pc{
(.,.)  & \relax : &    N\times M   \ar@{->}[rr]  && \bK [G]\\
 &&n,m   \ar@{|->}[rr] && (n,m) = \sum_{\sigma \in G} <m.\sigma^{-1},n>\sigma
}
  \end{equation}
  
\begin{proposition}
  The map $(.,.)$ in Equation~(\ref{eq:inf}) is $\bK[G]$-bilinear.
\end{proposition}
\proof Indeed
for any $\tau$ in $G$,  $m$ in $M$, and $n$ in $N$
\begin{eqnarray*}
  (\tau . n, m)&=&\sum_{\sigma \in G} <m.\sigma^{-1},\tau .n> \sigma\\
  &=& \sum_{\sigma \in G} <m.\sigma^{-1}\tau^{-1},\tau .n> \tau \sigma\\
  &=& \sum_{\sigma \in G} <m.\sigma^{-1}, n> \tau \sigma\\
  &=& \tau \sum_{\sigma \in G} <m.\sigma^{-1}, n>  \sigma\\
  &=& \tau (n,m).
\end{eqnarray*}

And

\begin{eqnarray*}
  (n, m . \tau )&=&\sum_{\sigma \in G} <m.\tau \sigma^{-1}, n> \sigma\\
  &=& \sum_{\sigma \in G} <m.\tau \tau^{-1}\sigma^{-1}, n>  \sigma\tau\\
  &=& \sum_{\sigma \in G} <m.\sigma^{-1}, n> \sigma \tau \\
  &=& (n,m)\tau.
\end{eqnarray*}

\hfill $\Box$

\subsection{Orthogonality}\label{sec:orth}

In the situation of Section~\ref{sec:inf} we consider
a right submodule $U$ of the  $\bK[G]$-module $M$. Call \[U^\perp
= \lbrace n\in N \,\,  | \,  <U,n>=0\rbrace \]
the orthogonal to $U$ in $N$ for the $<.,.>$ form. This is a $\bK$-vector space.
Since $U$ is stable by the action of $G$, its orthogonal $U^\perp$
is a left $\bK[G]$-module. And $U^\perp$ is the orthogonal
to $U$ for the $(.,.)$ form:
\[U^\perp
= \lbrace n\in N \,\, |  \,\, (n,U)=0\rbrace .\]
 We  consider similarly a  left $\bK[G]$-submodule $V$ of $N$
and let  \[V^\circ
= \lbrace m\in M \,\, | \,  <m,V>\, = 0\rbrace \]
be the orthogonal to $V$ in $M$ for the $<.,.>$ form. This
is a right $\bK[G]$-module. And $V^\circ$ is the orthogonal
to $V$ for the $(.,.)$ form:
\[V^\circ
= \lbrace m\in M \,\, | \,\,   (V,m)=0\rbrace . \]
We have $U\subset (U^\perp)^\circ$
and $V\subset (V^\circ )^\perp$. These
inclusions are equalities when  $M$ and $N$ are finite dimensional and 
 $<.,.>$   is perfect.

\subsection{The dual of a $\bK[G]$-module}\label{sec:duamo}
Let $N$ be a left  $\bK[G]$-module. We can see $N$ as a $\bK$-vector space and let  $\hat N$ be 
its  dual. This is 
a right  $\bK[G]$-module. For every
$\varphi$ in $\hat N$ and $\sigma$ in $G$
we set $\varphi.\sigma = \varphi \circ \sigma $. We consider the canonical $\bK$-bilinear form
defined by 
\[<\varphi , n > \, =\, \varphi (n)\]
for every $n$ in $N$ and $\varphi$ in $\hat N$. For every $\sigma $ in $G$
we have
\[<\varphi .\sigma , n> \, = \,  \varphi (\sigma . n) \, = \, <\varphi, \sigma .n>\]
so $<.,.>$ is invariant by $G$. Following Section~\ref{sec:inf} we define a $\bK[G]$-bilinear form
\begin{equation*}
  (.,.) : N\times \hat N \rightarrow \bK [G] \end{equation*}
by \begin{equation}\label{eq:isodu}(n,\varphi )=\sum_{\sigma \in G} \varphi(\sigma^{-1}.n)\sigma.\end{equation}

We define a  map from $\hat N$ to the dual $\check N$
of $N$ as a $\bK[G]$-module,
by sending $\varphi$ to the map \begin{equation}\label{eq:phi_G}\varphi^G : n\mapsto (n,\varphi).\end{equation}
We prove that this map is a bijection.
First  $\varphi \mapsto \varphi^G$
is trivially seen to be an injection.
As for surjectivity,  let $\psi  : N \rightarrow \bK[G]$ be a
$\bK[G]$-linear map. Writing \[\psi (n)=\sum_{\sigma \in G} \psi_\sigma (n)\sigma \]  we define
a $\bK$-linear coordinate form $\psi_\sigma$ on $N$ for every $\sigma $ in $G$.
We deduce from the $\bK[G]$-linearity of $\psi$
that $\psi_\sigma (n)=\psi_1(\sigma^{-1} . n)$ where  $1\in G$ is the identity element.
So  $\psi (n) = (n, \psi_1)$ for every $n$ in $N$. So $\psi = (\psi_1)^G$.


\subsection{Free submodules of a $\bK[G]$-module}\label{sec:freedi}

The ring $\bK [G]$ may not be semisimple. Still
free $\bK [G]$-submodules of finite rank are direct summands.

\begin{proposition}\label{prop:ds}
  Let $G$ be finite group,  $\bK$ a commutative field, and  $N$
  a left $\bK [G]$-module. Let   $V$  a  submodule
  of $N$. If $V$ is free   of finite rank  then
  it is a direct summand: there exists  a submodule $W$ of $N$ such that $N=V\oplus W$. Such  a $W$
  is called a complementary submodule to $V$.
\end{proposition}

\proof Let $r$ be the rank of $V$.
Let $v_1$, $v_2$, \ldots, $v_r$ be a basis of $V$. 
 Let $\varphi_1$, $\varphi_2$, \ldots,
$\varphi_r$ be the dual basis. For every $i$ such
that $1\leqslant i\leqslant n$, the coordinate form
$\varphi_{i,1}$ associated to the identity element $1$ in $G$
belongs to $\hat V$. Let $\psi_i$ be a $\bK$-linear form
on $N$ whose restriction to $V$ is $\varphi_{i,1}$.
Let $\psi_i^G \in \check N$ be the associated $\bK[G]$-linear form
according to Equations~(\ref{eq:isodu}) and (\ref{eq:phi_G}).
The restriction of  $\psi_i^G$ to $V$ is $\varphi_{i,1}^G$ and this  is $\varphi_i$.
The map
  \begin{equation*}
\xymatrix@R-2pc{
 \psi  & \relax : &    N   \ar@{->}[r]  & V\\
 &&n   \ar@{|->}[r] &
 \sum_{1\leqslant i\leqslant r}  \psi_i^G (n) .   v_i}
\end{equation*}is a $\bK[G]$-linear projection onto $V$.
Its kernel is a complementary $\bK[G]$-submodule to $V$.
\hfill $\Box$

\medskip

Proposition \ref{prop:ds} is a consequence of the Frobenius property which is known to be
satisfied by $\bK[G]$. See
\cite[Chapter IX]{currei}. The proof above provides an algorithm to compute a complementary module.

\section{Curves with a   group
  action}\label{sec:curveact}

Let $\bK$ be a  commutative field.
Let  $p$ be  the characteristic of $\bK$.
Let $X$ and $Y$
be two smooth, projective, absolutely integral  curves over $\bK$. Let $g_X$ be the genus of $X$ and let  $g_Y$ be the genus of $Y$.
Let $\tau : Y \rightarrow X$ be a Galois  cover with
Galois group $G$. Let $\co$ be the order of $G$.
There is a natural left action of $G$ on $\bK(Y)$ defined by
\begin{equation}\label{eq:actfunc} \sigma . f  = f\circ \sigma ^{-1} \text{\,\,\,\, for\,\,}
f\in \bK(Y) \text{\,\,  and  \,\,} \sigma \in G.\end{equation}  
There is a natural right action of $G$ on meromorphic
differentials  defined by
\begin{equation}\label{eq:actdif}\omega . \sigma = \sigma^\star \omega  \text{\,\,\,\, for\,\,}
\omega \in  \Omega_{\bK (Y)/\bK}   \text{\,\,  and  \,\,}
\sigma \in G.\end{equation} These are
$\bK(X)$-linear actions.
And the two actions
are compatible in the sense that 
\begin{equation}\label{eq:comp}
  (\omega . \sigma)(\sigma^{-1}. f)=(\omega f).\sigma
\end{equation}

We study  some free  $\bK[G]$-modules that arise
naturally in this context.

\subsection{The residue ring of a non-ramified fiber}\label{sec:rrf}

Let $P$ be a prime divisor (a place) on $X$.
Let $t_P$ be a uniformizing parameter at $P$.
Let \[a = \deg (P).\]
This is the degree  over $\bK$
of the residue field \[\bK_P=H^0(P,\cO_P)=H^0(X,\cO_X/\cO_X(-P)).\]
We assume that $\tau$ is not ramified above
$P$ and let $Q_1$ be a place above $P$.
Let  $G_1$ be
the decomposition group of $Q_1$. This is the stabilizer of $Q_1$
in $G$.
Places above $P$ are parametrized by left cosets  in $G/G_1$.
We write the fiber above $P$ \[Q=\sum_{\sigma \in G/G_1} Q_\sigma \text{\,\,\,\, with \,\,\,} Q_\sigma =\sigma (Q_1).\]
Let  \[b=[G:G_1]\] be the number
of places above $P$ and let \[c=\co /b=|G_1|\] be the residual degree, that is the degree
of \[\bK_\sigma =H^0(Q_\sigma, \cO_{Q_\sigma} )\] over   $\bK_P$ for all  $\sigma \in G/G_1$.
Let  \[\bR_Q = H^0(Q, \cO_Q)=H^0(Y, \cO_Y/\cO_Y(-Q))\]
be the residue ring  at $Q$. We have
\[\bR_Q = \bigoplus_{\sigma \in G/G_{1}}\, \bK_{\sigma }.\]
The action of  $G$ on $\bR_Q$
makes it  a free left $\bK [G]$-module of rank $a$.
Indeed it is a free $\bK_P [G]$-module  of rank $1$. A basis for it
consists of any normal element $\theta$ in 
 $\bK_{1}/\bK_P$.

If $m$ is a positive integer,  Taylor expansion
provides an isomorphism of  modules over $\bK_P [G]$
\[H^0(Y,\cO_Y/\cO_Y(-mQ)) \simeq \bR_Q[t_P]/t_P^m\]
between the residue ring at $mQ$ and the ring of truncated
series in $t_P$. So the former is a free left $\bK_P [G]$-module
of rank $m$.
A basis for it is made of the  $\theta t_P^k$ for $0\leqslant k < m$.

\subsection{The residue ring of   a non-ramified
  $G$-equivariant
  divisor}\label{sec:rrd}

We take $P$ an  effective divisor on $X$.
We assume that $\tau$ does not ramify above $P$
and let  $Q$ be the  pullback of $P$ by $\tau$. 
We write \[P=\sum_{1\leqslant i\leqslant I}m_iP_i.\]
Let $t_i$ be a uniformizing
parameter at $P_i$.
Let  $a_i$ be  the degree of the place $P_i$.
Let  $b_i$ be  the number of places of $Y$ above $P_i$.
Let $c_i=\co /b_i$. For every $1\leqslant i\leqslant I$
we choose a place $Q_{i,1}$ above $P_i$ and let  $G_{i,1}$ be 
the decomposition group at $Q_{i,1}$.
Let  $Q_i$ be  the pullback
of $P_i$ by $\tau$ and write
\begin{equation}\label{eq:actpoints}Q_i=\sum_{\sigma \in G/G_{i,1}}Q_{i,\sigma}  \text{\,\,\,\, with \,\,\,} Q_{i,\sigma}=\sigma (Q_{i,1})\end{equation} its decomposition
as a sum of $b_i$ places. Let $\bK_{i,\sigma}$ be  the residue field
at $Q_{i,\sigma}$.
We denote  by $\bA$ the residue algebra $H^0(Q,\cO_Q)$.
Taylor expansion induces an isomorphism of $\bK$-algebras
\begin{equation*}
  \bA=H^0(Q, \cO_Q)=H^0(Y, \cO_Y/\cO_Y(-Q))
  \simeq  \bigoplus_{i=1}^I\, \bigoplus_{\sigma \in G/G_{i,1}}\, \bK_{i,\sigma }[t_i]/t_i^{m_i}
  \end{equation*}
which is compatible
with the left actions of $G$ as defined by Equations (\ref{eq:actfunc}) and  (\ref{eq:actpoints}). 
In the special case when all the places $P_i$ have degree one,
a basis for 
$H^0(Q, \cO_Q)$ as
a  $\bK [G]$-module
is made of the $\theta_it_i^{k_i}$
for $1\leqslant i\leqslant I$  and $0\leqslant k_i < m_i$
where $\theta_i$ is  a normal element in the extension
$\bK_{i,1}/\bK$.
The proposition below follows
from the discussion in this section and the previous one.

\begin{proposition}\label{prop:freeres}
  Assume  the hypotheses at the beginning of Section~\ref{sec:curveact}.
  Let $P$ be an effective divisor
  on $X$. 
  Assume that $\tau$ is not ramified above $P$
  and let $Q$ be the pullback of $P$ by $\tau$.
  The residue ring
  $H^0(Q, \cO_Q)$
  is a free $\bK[G]$-module of rank the degree of $P$.
\end{proposition}  

\subsection{Duality}\label{sec:duaa}
We need a dual of
$\bA$ as a $\bK$-vector space.
We set \begin{equation*}\label{eq:hA}
  \hat \bA =  H^0(Y, \Omega_{Y/\bK}(-Q)/\Omega_{Y/\bK})   \simeq
  \bigoplus_{i=1}^I\, \bigoplus_{\sigma \in G/G_{i,1}}\,
  \left(   \bK_{i,\sigma }[t_i]/t_i^{m_i} \right) \frac{dt_i}{t_i^{m_i}}.\end{equation*}

For $f\in \bA$ and $\omega \in \hat\bA$ we write $<\omega,f>$
for the sum of the  residues of $\omega f$ at all the geometric points
of  $Q$.
This is a $\bK$-bilinear form.
We deduce from Equation~(\ref{eq:comp}) that this form is invariant by the action
of $G$ \[<\omega . \sigma , f > \,  = \,  <\omega, \sigma . f >\]
We define a $\bK[G]$-bilinear  form  using the
construction in Section~\ref{sec:inf}
\begin{equation}\label{eq:bil}
  (f,\omega)=\sum_{\sigma \in G}<\omega . \sigma^{-1}, f>\sigma \in \bK[G].\end{equation}

These two bilinear forms turn $\hat\bA$ into the dual
of $\bA$ as a $\bK$-vector space (resp. as a $\bK[G]$-module).
In the special case when all the places $P_i$ have degree one,
the dual basis to the basis introduced before Proposition \ref{prop:freeres}
is made of the $\mu_i t_i^{- k_i}dt_i/t_i$
for $1\leqslant i\leqslant I$  and $0\leqslant k_i < m_i$
where $\mu_i$ is the dual to the  normal element
$\theta_i$ 
in the extension
$\bK_{i,1}/\bK$.

\section{Free commutative actions}\label{sec:comm}

We study the situation at the beginning of Section~\ref{sec:curveact}
in the special case when the Galois cover $\tau : Y\rightarrow X$
is abelian and unramified. We prove that large enough equivariant
Riemann-Roch spaces
are free $\bK[G]$-modules. To this end we prove in Section \ref{sec:rrsp}
that
evaluation at some fibers induces  an isomorphism with one of the
 $\bK[G]$-modules
studied in Section \ref{sec:rrd}. We  need
a criterion for an  equivariant divisors on $Y$ to be non-special.
We recall  such  a criterion  in Section \ref{sec:sid}.
We introduce in Section \ref{sec:orths} the evaluation, interpolation
and  checking matrices whose existence follows
from the freeness of the considered modules.

\subsection{Special invariant divisors}\label{sec:sid}

The pullback by $\tau$ of a degree $g_X-1$ divisor
on $X$ is a degree $g_Y-1$ divisor on $Y$ according to the Riemann-Hurwitz formula.
We need a  criterion 
for the latter divisor
to be special. We will say that a divisor class is  {\it effective}
if it contains an effective divisor. When the degree of the class is the genus of the curve minus one, being
effective is equivalent to being special.

 \begin{proposition}\label{prop:special}
   Assume the
   hypotheses at the beginning of Section~\ref{sec:curveact}
   with $\tau$ abelian and unramified and $\bK$
   algebraically closed. Write $\co = \co_p\times \co_{p'}$
   where $\co_p$ is the largest power of $p$ dividing $\co$.
   Let    $c$ be a divisor class
   of degree $g_X-1$ on $X$
   and let  $\tau^\star(c)$ be  its pullback on $Y$. 
  If the  class $\tau^\star(c)$ is effective then $c$ is the sum of an effective class
  of degree $g_X-1$ 
  and a class of degree $0$ annihilated by $\tau^\star$ and by $\co_{p'}$.
\end{proposition}

\proof From \cite[\S 14]{coez}. Let $D$ be a divisor in $c$ and let $E$ be the pullback of $D$
by $\tau$. We assume that $\tau^\star(c)$ is effective. The space $H^0(Y, \cO_Y(E))$
is non-zero and is acted on by $G$.
Recall that a finite set of commuting endomorphisms
of a finite dimensional vector space over  an algebraically
closed field   has a common eigenvector. 
Let $f$ be such an eigenvector
for the  action of $G$. The divisor of $f$ is $J-E$ where $J$ is effective and stable
under the action  of $G$. So there exists an effective divisor $I$ on $X$ such
that $J$ is the pullback of $I$ by $\tau$. And the class of
$I-D$ is annihilated by $\tau^\star$. It is also annihilated by $\co_{p'}$ because
$f^{\co_{p'}}$ is invariant by $G$. \hfill $\Box$

\subsection{Riemann-Roch spaces}\label{sec:rrsp}

Let $E$ be a divisor on $Y$ defined over $\bK$ and invariant
by $G$.
The Riemann-Roch space $H^0(Y, \cO_Y(E))$ is a $\bK[G]$-module.
This module is free provided the degree of $E$ is large enough.

\begin{theorem}\label{prop:free}
Assume the hypotheses at
the beginning of Section~\ref{sec:curveact} with
$\tau$ abelian and unramified.
  Let $D$ be a divisor on $X$ with degree
  $\geqslant 2g_X-1$. Let $E$ be the  pullback of $D$
  by $\tau$.  
The $\bK$-vector space $H^0(Y, \cO_Y(E))$  is a free $\bK[G]$-module
  of rank $\deg (D)-g_X+1$.
\end{theorem}

\proof The statement is empty if
$g_X=0$. We assume that $g_X\geqslant 1$.
Because of
the Noether-Deuring theorem \cite[\S 2, Section 5]{bouralg8}, we  can
assume that $\bK$ is algebraically closed.
Let $k=\deg (D)-g_X+1$. We note that $k\geqslant g_X$. By dimension
count,  there exist
$k$ points 
\[P_1, P_2, \ldots, P_k \text{\,\,\, on \,\,} X\] 
such that the class
of $D-P_1-P_2-\dots -P_k$ is not the sum of an effective class
of degree $g_X-1$ and a class annihilated by
\[\tau^\star
  :\Pic (X)\rightarrow \Pic(Y).\]
Indeed every divisor class of degree $g_X-1$ contains a divisor
of the form $D-P_1-P_2-\dots -P_k$ because $k$ is greater  than
or equal to the dimension $g_X$ of $\Pic^{g_X-1}$. On the other
hand the set of effective classes of degree $g_X-1$ has dimension
$g_X-1$. And the kernel of $\tau^*$ is finite. So the set of 
bad classes has codimension $1$ in $\Pic^{g_X-1}$.

Let $P$ be the divisor 
sum of all $P_i$ and let $Q$ be its pullback  by $\tau$.
According to
Proposition \ref{prop:special} the class of $E-Q$ is ineffective. Thus the evaluation
map \[H^0(Y, \cO_Y(E))\rightarrow H^0(Y, \cO_Y(E)/\cO_Y(E-Q))\simeq H^0(Q, \cO_Q)\] is an isomorphism of $\bK[G]$-modules.
Proposition~\ref{prop:freeres} then implies
that $H^0(Q, \cO_Q)$ is a free $\bK[G]$-module of rank $k$. \hfill $\Box$

\medskip

When the degree of $D$ is smaller than $2g_X-1$ it is not granted that
$H^0(Y, \cO_Y(E))$ is free. We mention two useful  partial results.

\begin{proposition}\label{prop:p'}
Assume the hypotheses at
the beginning of Section~\ref{sec:curveact}  with
$\tau$ abelian and unramified. Assume that $p$ does not divide
$\co$.
  Let $D$ be a divisor on $X$ with degree
  $\geqslant g_X$. Let $E$ be the  pullback of $D$
  by $\tau$.  
Then $H^0(Y, \cO_Y(E))$  contains  a free $\bK[G]$-module
of rank $\deg (D)-g_X+1$.
\end{proposition}
\proof The ring $\bK[G]$ is semi-simple.
Let $\cL(E)=H^0(Y, \cO_Y(E))$.
Let $m$ be the smallest among the multiplicities in $\cL(E)$
of irreducible
representations of $G$. This is
the smallest among the multiplicities of multiplicative
characters of $G$ in $\cL(E)\otimes \bKb$
where $\bKb$ is an algebraic closure of $\bK$.
It is clear that $\cL(E)$ contains $m$ copies
of the regular representation of $G$. On the other hand
let $\chi : G\rightarrow \bKb$ be a multiplicative character.
By the normal basis theorem
there exists  an eigenfunction $r$  in $\bKb(Y)$ associated with $\chi$.
The divisor of $r$  is the pullback by $\tau$ of a divisor
$R$ on $X$.
Let $\cL(E)_\chi$ be the eigenspace in $\cL(E)$ associated with $\chi$.
The map $f\mapsto f/r$ is a bijection between $\cL(E)_\chi$
and $H^0(X, \cO_X(D+R))$. The dimension of the latter is at least 
$\deg (D)-g_X+1$. Thus  $m\geqslant \deg (D)-g_X+1$.
\hfill $\Box$

\medskip

We can say something also when $G$ is a $p$-group and  $\bK$  a finite field.

\begin{proposition}\label{prop:pfin}
Assume the hypotheses at
the beginning of Section~\ref{sec:curveact}  with
$\tau$ abelian and unramified. Assume that $\bK$ is a finite field
with at least four elements. Assume that $\co$ is a power of $p$.
Assume that $g_X\geqslant 2$.
Let $d\geqslant g_X$ be an integer.
Let $r=d-g_X+1$. Assume that 
there exists an effective divisor on $X$
with degree $r$ and defined over $\bK$.
Then there exists a divisor  $D$  on $X$
such that $D$ is defined
over $\bK$, $D$ has  degree $d$, and
$H^0(Y, \cO_Y(E))$  is  a free $\bK[G]$-module
of rank $r=d-g_X+1$ where $E$ is  the  pullback of $D$
  by $\tau$.  
\end{proposition}  
  
\proof Set $r=d-g_X+1$. Let $P$ be an effective divisor on $X$
with degree $r$ and defined over $\bK$. According
to \cite[Theorem 11]{BLB} by  Ballet and Le Brigand,
there exists a degree $g_X-1$
non-special divisor $I$ defined over $\bK$.
Set $D=I+P$. Let $E$, $J$, and $Q$ be the pullbacks
of $D$, $I$, and $P$ by $\tau$.
The class of the divisor $J$ is ineffective according to Proposition \ref{prop:special}. So the
evaluation map $H^0(Y, \cO_Y(E))\rightarrow H^0(Q, \cO_Q)$ is a bijection.
And the latter is a free $\bK[G]$-module
according to Proposition~\ref{prop:freeres}. \hfill $\Box$

\medskip

Theorem \ref{prop:free} translates into  a similar statement for differentials.

\begin{proposition}\label{prop:freediff}
Assume the hypotheses at
the beginning of Section~\ref{sec:curveact} with
$\tau$ abelian and unramified.
  Let $D$ be a divisor on $X$ with $\deg D <0$. Let $E$ be the  pullback of $D$
  by $\tau$.  
The $\bK$-vector space $H^0(Y, \Omega_{Y/\bK}(E))$  is a free $\bK[G]$-module
  of rank $g_X-1-\deg (D)$.
\end{proposition}

\proof
The statement is trivial if $g_X=0$. We assume that $g_X\geqslant 1$.
Let $\omega_0$ be a non-zero holomorphic differential on $X$.
The pullback of $\omega_0$ on $Y$ by $\tau$ is denoted by
$\omega_0$ also. The map $\omega \mapsto \omega/\omega_0$ is an isomorphism of
$\bK$-vector spaces between $H^0(Y, \Omega_{Y/\bK}(E))$ and $H^0(Y, \cO_Y( (\omega_0)-E))$. According to Equation
(\ref{eq:comp}) this isomorphism  is compatible
with the actions of $G$ on either sides given by Equations (\ref{eq:actfunc}) and (\ref{eq:actdif}).
Since the degree of $(\omega_0)-D$ is at least $2g_X-1$  we can apply Theorem \ref{prop:free} to prove 
that $H^0(Y, \cO_Y( (\omega_0)-E))$ is free and deduce that $H^0(Y, \Omega_{Y/\bK}(E))$ is free as well. \hfill $\Box$

\subsection{The orthogonal submodule}\label{sec:orths}

In the situation of the beginning of
Section~\ref{sec:curveact} and assuming that 
$\tau$ is abelian and unramified  we let $D$ and $P$ be  divisors on $X$
with  $P$ effective.
We assume that $D$ and $P$ are disjoint.
We assume that \begin{equation}\label{eq:inte}2g_X-1\leqslant \deg (D) \leqslant \deg (P) -1.\end{equation}
Let  $E$ be  the pullback of $D$
by $\tau$ and  let  $Q$ be  the pullback of $P$. 
We write
\[\cL (E)=H^0(Y, \cO_Y (E))  \text{\,\,\, and \,\,\,}  \itOmega (-Q+E)=
H^0(Y, \Omega_{Y/\bK}(-Q+E)).\]
Theorem \ref{prop:free}, Proposition \ref{prop:freediff},  and Equation (\ref{eq:inte}) imply that these two $\bK[G]$-modules 
are free. And 
the evaluation  maps 
\[\cL (E)\longrightarrow \bA \text{\,\,\, and \,\,\,} \itOmega (-Q+E) \longrightarrow \hat\bA \text{\,\,\, are injective.}
\]
So $\cL(E)$ can be seen as
 a free submodule of $\bA$ and
 $\itOmega (-Q+E)$ as
  a free submodule of $\hat\bA$.
For dimension reasons and
 due to the residue theorem, these two $\bK[G]$-modules are orthogonal to each other for the form introduced in Equation~(\ref{eq:bil}).
Proposition \ref{prop:ds} implies that $\cL(E)$  has a complementary
submodule in $\bA$ that is isomorphic to the dual of  $\itOmega (-Q+E)$ and is thus a free submodule.
Similarly  $\itOmega (-Q+E)$ has a free complementary submodule in $\hat\bA$ that is
isomorphic to  the dual of $\cL(E)$.

In the special case when all the places $P_i$ have degree one,
we have introduced a
natural basis for $\bA$
before Proposition \ref{prop:freeres}
and its dual basis $\hat\bA$
in Section \ref{sec:duaa},
using Taylor expansions at the places above the $P_i$.

We choose $\bK[G]$-bases for $\cL (E)$ and  $\itOmega (-Q+E)$.
We denote by $\cE_E$ the $\deg (P)\times (\deg (D)-g_X+1)$
matrix with coefficients in $\bK [G]$ of  the evaluation map $\cL (E)\rightarrow \bA$
in the chosen bases.
We   denote by $\cC_E$ the $\deg (P)\times (\deg (P)-\deg (D)+g_X-1)$
matrix 
of  the  map $\itOmega (-Q+E)\rightarrow \hat\bA$
in the chosen  bases.
The matrix  $\cC_E$ checks that a vector in $\bA$ belongs to $\cL(E)$. Its left kernel is the image of $\cE_E$.
So \[\cC_E^t\times \cE_E = 0,\]
a zero $(\deg (P)-\deg (D)+g_X-1)\times (\deg (D)-g_X+1)$
matrix with entries in $\bK[G]$. 

\smallskip
We choose a $\bK[G]$-linear
projection $\bA\rightarrow \cL (E)$ and denote by $\cI_E$
the  matrix  of this projection.  This is a $(\deg (D)-g_X+1)\times \deg (P)$
matrix with coefficients in $\bK[G]$.
This is an interpolation matrix since it recovers a function in $\cL (E)$
from its evaluation at $Q$. Equivalently
\[\cI_E\times \cE_E=1\]
the $(\deg (D)-g_X+1) \times (\deg (D)-g_X+1)$ identity matrix
with coefficients in $\bK[G]$.
We note that applying either of the matrices $\cE_E$, $\cC_E$, $\cI_E$ requires at most a constant times
$\deg (P)^2$  operations in $\bK[G]$.

\section{Pad\'e approximants}\label{sec:pade}

In the situation of the beginning of
Section~\ref{sec:curveact} and assuming that 
$\tau$ is abelian and unramified
we let $D_0$, $D_1$  and $P$ be  divisors on $X$
with  $P$ effective.
We assume that $D_0$ and $D_1$ are  disjoint from $P$.
Let  $E_0$, $E_1$, and $Q$  be  the pullbacks of $D_0$, $D_1$,
and $P$ 
by $\tau$.
We assume that  
\begin{equation}\label{eq:pade1}
  2g_X-1\leqslant  \deg (D_1) \leqslant \deg (P) -1,\end{equation}
\begin{equation}\label{eq:pade2}
  g_X\leqslant \deg (D_0)\leqslant \deg (P) -1.\end{equation}
Equation (\ref{eq:pade1}) implies that
the $\bK[G]$-modules  $\cL (E_1)$ and  $\itOmega (-Q+E_1)$ are free and 
the evaluation  maps into $\bA$ and  $\hat\bA$  are injective.
We assume that $\cL(E_0)$ contains a free $\bK[G]$-module
of rank $\deg (D_0)-g_X+1$ and denote by $\cL(E_0)_\bff$ such a submodule.

\smallskip
Given  $r$ in $\bA$, $a_0\not = 0$ in $\cL (E_0)$ and $a_1$
in $\cL (E_1)$
such  that
\[a_0r-a_1=0 \in \bA,\]
we say that $(a_0, a_1)$ is a {\bf Pad\'e approximant} of $r$
and  we say that   $a_0$ is   a {\bf denominator} for $r$.
Denominators for $r$ are non-zero $a_0$ in $\cL (E_0)\subset \bA$ such that
\[a_0r\in \cL (E_1).\] Equivalently
\begin{equation}\label{eq:coliK}
  (a_0r,\omega)=0 \text{\,\,\, for every \,\,\,}
\omega \in \itOmega (-Q+E_1).\end{equation}  Denominators are thus non-zero solutions
of a $\bK$-linear system of equations. We note that this is  not a $\bK[G]$-linear system in general because
multiplication by $r$ is not $\bK[G]$-linear. In
Section \ref{sec:splitc} we show that one can be a bit more explicit
in some cases. We consider the problem of computing Pad{\'e} approximants
in Section \ref{sec:wied}.

\subsection{The split case}\label{sec:splitc}

Assume  that $P=P_1+\dots+P_n$
is a sum of $n$ pairwise distinct rational  points over $\bK$. Assume that the fiber of $\tau$
above each $P_i$ decomposes as a sum of $\co$ rational points over $\bK$. We choose a point $Q_{i,1}$
above each $P_i$ and set
\[Q_{i,\sigma }= \sigma (Q_{i,1}) \text{\,\,\, for every \,\, } \sigma \in G.\]
  For every $1\leqslant i\leqslant n$  let  $\alpha_i$ be  the function in $\bA$ that takes value
  $1$ at $Q_{i,1}$ and zero everywhere else. We thus form a
  basis \[\cA_G= (\alpha_i)_{1\leqslant i
    \leqslant n}\] of $\bA$ over $\bK[G]$.
  We note $\hat\cA_G$ its dual basis.
For every $1\leqslant i\leqslant n$ and $\sigma \in G$
let
\[\alpha_{i,\sigma } =\sigma . \alpha_i = \alpha_i \circ \sigma^{-1}\] be the function in $\bA$ that takes value
$1$ at $Q_{i,\sigma }$ and zero  everywhere else.
We thus form a basis \[\cA_\bK=(\alpha_{i,\sigma })_{1\leqslant i\leqslant n, \, \sigma  \in G}\] of $\bA$ over $\bK$.
The coordinates 
of $r$ in the $\bK[G]$-basis $\cA_G$ 
are \[r_G=(\sum_{\sigma \in G }r(Q_{i,\sigma})\sigma )_{1\leqslant i \leqslant n}\] and the coordinates 
of $r\in \bA$ in the $\bK$-basis $\cA_\bK$ 
are \[r_\bK=(r(Q_{i,\sigma }))_{1\leqslant i\leqslant n, \, \sigma \in G}.\]

Multiplication by $r$ is a $\bK$-linear map
from $\bA$ to $\bA$. Let  \[\cR_\bK \in \cM_{\co . n, \co . n}(\bK)\]
be the $\co . n\times \co . n$ diagonal matrix
of this map in the basis $\cA_\bK$.

We choose a $\bK[G]$-basis $\cZ_G$ for $\cL (E_0)_\bff$ and denote by $\cE_G^0$ the
$\deg (P)\times (\deg (D_0)-g_X+1)$ matrix of the
$\bK[G]$-linear injective map \begin{equation}\label{eq:cL}
  \cL (E_0)_\bff\rightarrow \bA\end{equation} in the  bases
$\cZ_G$ and $\cA_G$.
We denote by $\cZ_\bK$ the $\bK$-basis of $\cL (E_0)_\bff$ obtained
by letting $G$ act on $\cZ_G$.
Let   $\cE_\bK^0$ be the
matrix of the map (\ref{eq:cL}) 
in the  bases
$\cZ_\bK$ and $\cA_\bK$.
The matrix
$\cE_\bK^0$ is obtained from
$\cE_G^0$ by replacing each $\bK[G]$ entry by the corresponding
$\co\times \co$ circulant-like matrix with entries in  $\bK$.

Let
$\hat \cA_G$ be the basis of the $\bK[G]$-module $\hat\bA$,
dual to $\cA_G$.
We choose a $\bK[G]$-basis $\cU_G$
for $\itOmega (-Q+E_1)$ and denote by $\cC_G^1$
the   matrix of the
injective map \begin{equation}\label{eq:omeg}
  \itOmega (-Q+E_1)\rightarrow \hat\bA\end{equation} in the  bases $\cU_G$
and $\hat \cA_G$.
This is a 
$\deg (P)\times (\deg (P)-\deg (D_1)+g_X-1)$  matrix with entries in $\bK[G]$.
Let  $\cU_\bK$ be  the $\bK$-basis of $\itOmega (-Q+E_1)$
obtained  by letting $G$ act on  $\cU_G$.
Let
$\hat \cA_\bK$ be the basis of the $\bK$-vector space
$\hat\bA$,
dual to $\cA_\bK$.
The 
matrix of 
the map (\ref{eq:omeg})
in the  bases
$\cU_\bK$ and $\hat\cA_\bK$ is called $\cC_\bK^1$.

\smallskip

Let $a_0$ in  $\cL (E_0)_\bff$ and
let $x_G$ be the  coordinates of $a_0$
in the $\bK[G]$-basis $\cZ_G$. This is a column
of height $\deg(D_0)-g_X+1$. We let  $x_\bK$
be the coordinates of $a_0$
in the $\bK$-basis $\cZ_\bK$. This is a column
of height $\co . (\deg(D_0)-g_X+1)$ obtained from
$x_G$ by replacing each entry by its  $\co$ coefficients
in the canonical basis of $\bK[G]$.
We deduce from Equation (\ref{eq:coliK}) that $a_0$
is a denominator for $r$
if and only if $x_\bK$
is in the kernel of the matrix \[\cD_r=(\cC_\bK^1)^t\times \cR_\bK\times \cE^0_\bK
\in \cM_{\co .  (\deg P -deg D_1+g_X-1)\times \co . (\deg D_0-g_X+1)}(\bK).\]

\begin{proposition}\label{prop:testde}
  Assume that we are in  the context of the beginning of Section \ref{sec:pade}. In particular assume  Equations
  (\ref{eq:pade1}) and   (\ref{eq:pade2}),  assume
  that $P$ is a sum of $n$ pairwise
  distinct $\bK$-rational points,   and that  the $n$ corresponding fibers of $\tau$ split over $\bK$.
  Assume that we are given the matrices $\cE^0_\bK$ and $\cC_\bK^1$.
  On input an   $r = (r(Q_{i,\sigma }))_{1\leqslant i\leqslant n, \, \sigma  \in G}$
  in $\bA$
  and some
  $a_0$ in $\cL(E_0)_\bff$, given by its coordinates $x_\bK$ in the basis
  $\cZ_\bK$,  one can check if  $a_0r\in  \cL(E_1)$
  at the expense of $\cQ .  n^2$ operations  in $\bK[G]$ (addition, multiplication) and $\cQ . \co. n$ operations in $\bK$ (addition, multiplication) where $\cQ$ is some absolute constant.
\end{proposition}

\proof We first multiply $x_\bK$
by $\cE^0_\bK$ or rather $x_G$ by $\cE^0_G$. This 
requires less than $2\deg (P)\times (\deg (D_0)-g_X+1)$ operations in $\bK[G]$.
We then multiply the result by $\cR_\bK$. This
requires less than $\co . \deg (P)$ operations in $\bK$ because $\cR_\bK$ is diagonal.
We finally multiply the result by $(\cC_\bK^1)^t$ or rather $(\cC_G^1)^t$.
This  requires less than $2\deg (P)\times (\deg (P)-\deg (D_1)+g_X-1)$ operations  in $\bK[G]$. \hfill $\Box$

\subsection{Computing Pad{\'e} approximants}\label{sec:wied}

According to Proposition \ref{prop:testde} one can  efficiently check a denominator.
As a consequence, one can  find  a random
denominator, assuming  that there is some  in $\cL(E_0)_\bff$, using an iterative method as in
\cite{Wied, Kal}.
Recall that an $\ell\times n$ {\bf black box} matrix $A$ with coefficients
in a field $\bK$ is an oracle that on input
 an $n\times 1$ vector $x$ returns $Ax$.

\begin{proposition}[Wiedemann, Kaltofen, Saunders]\label{prop:wks}
  There exists a probabilistic (Las Vegas) algorithm
  that takes as input an $\ell\times n$  black box matrix $A$
  and an $\ell \times 1$ vector  $b$ with entries in
  a finite field $\bK$ and returns a uniformly distributed
  random solution $x$ to the system $Ax=b$, if there is some, with probability
  of success $\geqslant 1/2$ at the expense of $\cQ .  m . \log m$
  calls to the black box for $A$ and $\cQ.  m^2 . (\log (m))^2$
  operations in $\bK$ (addition, multiplication, inversion, picking a random element) where $\cQ$ is some absolute constant and
  $m=\max (\ell, n)$.
\end{proposition}

Using  Proposition \ref{prop:wks} and bounding the cost of a call to the black box with the help of Proposition
\ref{prop:testde}  we deduce 

\begin{proposition}\label{prop:finddenom}
  Under the hypotheses of Proposition \ref{prop:testde}
  and
  on input a vector $r = (r(Q_{i,j})_{i,j})$ in $\bA$
  one can find
  a uniformly distributed random
  denominator for $r$, if there is some in $\cL(E_0)_\bff$,
  with probability of success $\geqslant 1/2$, at
  the expense of $\cQ .  \co . n^3. \log (\co . n)$ operations
  in $\bK[G]$ (addition, multiplication) and $\cQ .  (\co . n .
  \log (\co . n))^2$ operations
  in $\bK$ (addition, multiplication, inversion, picking a random element)
  where $\cQ$ is some absolute constant.
\end{proposition}

Once we have found a denominator  $a_0$ for $r$
we set $a_1=ra_0$ and recover the coordinates of $a_1$
applying the interpolation matrix associated to $E_1$.

\section{Computing in the group algebra}\label{sec:ga}

Given a finite commutative group $G$ and a finite  field
$\bK$ we will need efficient algorithms to multiply
in $\bK[G]$. This is classically achieved using  a discrete
Fourier transform when $G$ is cyclic and $\bK$ contains
enough roots of unity. The complexity
analysis requires some care in general. This is the purpose of this section.
We recall in Section \ref{sec:ft} the definition of the  Fourier transform in the setting of commutative finite groups.
The most classical case of cyclic groups is studied in Section \ref{sec:unift} from an algorithmic
point of view.
The general case follows by induction as explained in Section \ref{sec:multiv}. The complexity of the resulting multiplication
algorithm in $\bK[G]$ is bounded in Section \ref{sec:fasmul}.
\subsection{Fourier transforms}\label{sec:ft}

Let $G$ be a finite commutative group.
Let $\co$ be the order of $G$. Let $e$
be its exponent.
Let $\bK$ be a commutative field  containing  a primitive
$e$-th root of unity. In particular
$e$ and $\co$ are non-zero  in  $\bK$.
Let  $\hat G$
be the dual of $G$ defined as the group of characters
$\chi : G\rightarrow \bK^*$.
We define a  map from
the group algebra of $G$ to the algebra of functions on $G$
  \begin{equation*}
\xymatrix@R-2pc{
 \top  & \relax : &    \bK[G]   \ar@{->}[r]  & \HomS (G, \bK)\\
 &&\sum_{\sigma \in G} a_\sigma \sigma  \ar@{|->}[r] & \sigma \mapsto a_\sigma }
  \end{equation*}
  This is an isomorphism of $\bK$-vector spaces. We let $\bot :
  \HomS (G, \bK)\rightarrow \bK[G]$
be   the reciprocal map.
We dually define 
  \begin{equation*}
\xymatrix@R-2pc{
 \hat \top  & \relax : &    \bK[\hat G]   \ar@{->}[r]  & \HomS (\hat G, \bK)\\
 &&\sum_{\chi \in \hat G} \, a_\chi \chi  \ar@{|->}[r] & \chi \mapsto a_\chi }
  \end{equation*}
  and its reciprocal map $\hat \bot$.
  We let  \[\iota_G : \bK[G]\rightarrow \bK[G]\] be the
  $\bK$-linear involution that maps $\sigma$
  onto $\sigma^{-1}$.
We define the Fourier transform 
  \begin{equation*}
\xymatrix@R-2pc{
 \FT_G    & \relax : &    \bK[ G]   \ar@{->}[r]  & \HomS (\hat G, \bK)\\
 &&\sum_{\sigma \in G} a_\sigma \sigma  \ar@{|->}[r] & \chi \mapsto \sum_{\sigma}a_\sigma \chi (\sigma) }
  \end{equation*}
  The Fourier transform evaluates an element  in the group algebra
  at every character. The Fourier transform of the dual group
  \begin{equation*}
\xymatrix@R-2pc{
 \FT_{\hat G}    & \relax : &    \bK[ \hat G]   \ar@{->}[r]  & \HomS (G, \bK)\\
 &&\sum_{\chi \in \hat G} a_\chi \chi  \ar@{|->}[r] & \sigma \mapsto \sum_{\chi}a_\chi \chi (\sigma) }
  \end{equation*}
  provides an inverse for $\FT_G$ in the sense that
  \[\bot\circ \FT_{\hat G}\circ \hat\bot \circ \FT_G = \co . \iota_G\]
  is the $\bK$-linear invertible map that sends  $\sigma$ to $\co . \sigma^{-1}$.

  Let $M$ be a finite dimensional $\bK$-vector space. We
  set \[M[G]=M\otimes_\bK \bK[G]\] and note that
  \[\HomS (\hat G,M)=M\otimes_\bK \HomS(\hat G, \bK).\]
We define a  Fourier transform on $M$
   \begin{equation*}
\xymatrix@R-2pc{
 \FT_M    & \relax : &     M [ G]   \ar@{->}[r]  & \HomS (\hat G, M)\\
 &&\sum_{\sigma \in G} m_\sigma \otimes \sigma  \ar@{|->}[r] & \chi \mapsto \sum_{\sigma} \chi (\sigma) m_\sigma  }
\end{equation*}
It turns a free $\bK[G]$-module into a free $\HomS (\hat G, \bK)$-module.

  \subsection{Univariate Fourier transforms}\label{sec:unift}

  We assume in this section that the group $G$ is cyclic
  of order $\co$. We choose a primitive
  $\co$-th root of unity $\omega$ in $\bK$. We choose a generator
  in $\hat G$ and deduce the following identifications
 \[\HomS (\hat G, \bK)=\bK^\co  \text{\,\,\, and \,\,\,} \bK[G]=\bK [x]/(x^\co-1).\]
  Let $M$ be a finite dimensional  $\bK$-vector space.
  Setting \[M[x]=M\otimes_\bK \bK[x] \text{\,\,\, and \,\,\,} M[G]=M\otimes_\bK \bK [x]/(x^\co-1).\]
   the Fourier transform is
\begin{equation*}
\xymatrix@R-2pc{
\FT_M   & \relax : &    M[G]  \ar@{->}[r]  & M^\co\\
&  & m  \ar@{|->}[r] & (m(1),m(\omega), m(\omega^2), \ldots, m(\omega^{\co-1}))}
\end{equation*}
Given $m$ in $M[G]=M\otimes_\bK \bK[x]/(x^\co-1)$ the computation of $\FT_M(m)$ reduces to the multiplication of a polynomial of degree $2\co-2$ in $\bK[x]$
and a vector  of degree $\co-1$ in $M[x]$ using formulae by
Rabiner, Schafer, Rader, and  Bluestein \cite{RSR,Blu}. 

\begin{proposition}\label{prop:ffcy}
  Let $\bK$  be a commutative field.
  Let $M$ be a  finite dimensional $\bK$-vector space.
Let $\co\geqslant 2$ be an integer.
Assume that $\bK$ contains a primitive
$\co$-th root of unity $\omega$
and a primitive root of unity of order a power
of two that is bigger than  $3\co-3$.
Let \[m=m_0\otimes 1 +m_1 \otimes x+\dots +m_{\co-1}\otimes x^{\co-1}\bmod x^\co-1 \in M\otimes_\bK \bK[x]/(x^\co-1).\]
One can compute $\FT_M(m)$ at the expense of
$\cQ . \co . \log \co$ additions, multiplications and
inversions in $\bK$, additions
and scalar multiplications in $M$, where $\cQ$
is an absolute constant.
\end{proposition}
\proof We adapt
the notation  from  \cite[I.5.4, Proposition 5.10]{bost}.
For every $0\leqslant i \leqslant 2\co-2$
let \[t_i=i(i-1)/2 \text{ \,\,\,  and  \,\,\, } \beta_i=\omega^{t_i}.\]
We note
that \[t_{i+1}=t_i+i \text{\,\,\,  and \,\,\, }\beta_{i+1}=\beta_i\omega^i.\]
So
one can compute the $\beta_i$  for
$0\leqslant i\leqslant 2\co -2$
at the expense of $4\co$
operations in $\bK$.
We then compute the inverse of every $\beta_i$.
For every $0\leqslant i \leqslant \co-1$ let
\[n_i=\beta_i^{-1}m_i.\] These can be computed at the
expense of $\co$ scalar multiplications in $M$.
Let \[n(x)=n_{\co -1}+n_{\co -2}\otimes x +\dots+n_{0}\otimes x^{\co -1} \in M[x]\]
and let  \[b(x)=\beta_0+\beta_1x+\dots +\beta_{2\co -2}x^{2\co -2} \in \bK[x].\]
Let \[r(x)=b(x). n(x) = \sum_{0\leqslant i\leqslant 3\co -3}
r_i \otimes x^i \in M[x].\]
From the identity \[t_{i+j}=t_i+t_j+ij\]
we deduce \[\omega^{ij}\beta_i\beta_j=\beta_{i+j} \text{\,\, for \,\,}  0\leqslant i, j \leqslant \co -1\]
and 
\[\sum_{j=0}^{\co -1}\omega^{ij}m_j=
\beta_i^{-1}\sum_{j=0}^{o-1}\beta_{i+j}n_j.\]
We deduce
 that 
$\FT_M(m)=(\beta_0^{-1}r_{\co -1}, \beta_1^{-1}r_{\co }, \beta_2^{-1}r_{\co +1}, \dots,
\beta_{\co -1}^{-1}r_{2\co -2})$.
Since 
$\bK$ contains a primitive root of unity of order a power
of two that is bigger than  $3\co -3$,  the coefficients
in the product
$r(x)=b(x).n(x)$
can be computed at the expense of $\cQ.\co.\log \co$
operations in $\bK$, additions in $M$ and products of a vector  in
$M$ by a scalar in $\bK$. See \cite[I.2.4, Algorithme 2.3]{bost}.
\hfill $\Box$

\subsection{Multivariate Fourier transforms}\label{sec:multiv}

Let $(\co _i)_{1\leqslant i\leqslant I}$ be integers such
that $2\leqslant \co _1  | \co_2 | \dots | \co_I$. Let $C_i=\bZ/\co_i\bZ$
and $G=\Pi_{1\leqslant i\leqslant I}C_i$. For $1\leqslant i\leqslant I$ set
\[A_i=\bK[C_i] \text{\,\,\, and \,\,\,} B_i=\HomS(\hat C_i,\bK).\]
For $0\leqslant i\leqslant I$ set
 \[M_i=\bigotimes_{j\leqslant i}B_j\otimes
  \bigotimes_{j> i}A_j.\]
So $M_0=\bK[G]$ and $M_I=\HomS(\hat G,\bK)$. For $0\leqslant i\leqslant I-1$
write \[M_i=\bigotimes_{j\leqslant i}B_j\otimes \bK[C_{i+1}]\otimes   \bigotimes_{j> i+1}A_j\]
as a $\bK[C_{i+1}]$-module and let  \[F_i : M_i\rightarrow M_{i+1}\]
be the corresponding Fourier transform as defined in Section
\ref{sec:unift}. We  check that \[\FT_G = F_{I-1}\circ F_{I-2}\circ
\dots \circ F_0.\]

Using Proposition \ref{prop:ffcy} we deduce

\begin{proposition}\label{prop:ffG}
Let $(\co_i)_{1\leqslant i\leqslant I}$ be integers such
that $2\leqslant \co_1  | \co_2 | \dots | \co_I$. Let
$G=\prod_{1\leqslant i\leqslant I}(\bZ/\co_i\bZ)$.
Let $\co$ be the order of $G$.
Let $e=\co_I$ be the exponent of $G$.
Let $\bK$ be a commutative field containing a primitive
root of unity of order $e$
and a primitive root of unity of order a power
of two that is bigger than  $3e-3$.
Given an element $a = \sum_{\sigma \in G}a_\sigma \sigma$ in $\bK[G]$
one can compute $\FT_G(a)$ in $\HomS (\hat G, \bK)$ at the expense
of $\cQ.\co.\log \co$ additions, multiplications and
inversions in $\bK$. Here $\cQ$ is some absolute constant.
\end{proposition}

\subsection{Fast multiplication in $\bK[G]$}\label{sec:fasmul}

Let $G$, $\co$, $e$ be as in Section \ref{sec:multiv}. 
Let $\bK$ be a commutative field.
In this section we study
the algorithmic complexity  of computing the product of two
 given elements
\begin{equation}\label{eq:ab}a =\sum_{\sigma \in G}a_\sigma  \sigma\text{ \,\,\, and \,\,   } b
  =\sum_{\sigma \in G}b_\sigma  \sigma  \text{\,\, in \,\, }  \bK[G].\end{equation}
It will depend on the field $\bK$. We first treat the case
  when $\bK$ has enough roots of unity.

\begin{proposition}\label{prop:enro}
  In the context of the beginning of Section \ref{sec:fasmul}
  assume that $\bK$ contains a primitive
root of unity of order $e$
and a primitive root of unity of order a power
of two that is bigger than  $3e-3$. One can compute the product
$ab\in \bK[G]$ at the expense of $\cQ.\co.\log \co$ operations in $\bK$
where $\cQ$ is some absolute constant.
\end{proposition}

\proof We compute $A=\FT_G(a)$ and $B=\FT_G(b)$ as in Section \ref{sec:multiv}.
We then compute  $C=AB$ in $\HomS (\hat G, \bK^*)$
at the expense of $\co$ multiplications in $\bK$. We then deduce
$c=ab$ applying $\FT_G^{-1}$ to $C$. The cost of this computation
is bounded
using Proposition \ref{prop:ffG}.\hfill $\Box$

\smallskip

We now consider the case when $\bK$ is $\bZ/p\bZ$
where $p$ is a prime integer. We miss roots of unity
in $\bK$ in general. So we transport the problem
into another ring using non-algebraic maps.
We let $t$ be the smallest
power of $2$ that is bigger than $3e-3$. Let $p'$
be the smallest prime integer congruent to $1$
modulo $\co . (p-1)^2.t$. We set $\bK'=\bZ/p'\bZ$ and note
that $\bK'$ contains a primitive root of unity of order
$e$ and a primitive root of order a power of two bigger
than $3e-3$. Also \begin{equation*}\label{eq:pq} p' > \co . (p-1)^2.\end{equation*}
By a result of Heath-Brown,  the exponent in  Linnik's theorem for primes
in arithmetic progressions can be taken to be $11/2$. See
\cite{hb} and the  recent improvement \cite{xyl}.
We deduce that there exists an
absolute constant $\cQ$ such that
\begin{equation*}\label{eq:linn} p' \leqslant  \cQ
  (\co . p)^{11}.\end{equation*} 
For $c$ a congruence class in $\bK=\bZ/p\bZ$ we
denote by  $\ell (c)$ the lift of $c$, that is   the unique integer
in the intersection
of  $c$ with the interval $[0,p[$. We write \begin{equation}\label{eq:up}\up (c)=\ell (c)\bmod p'.\end{equation}
We thus define   maps $\ell : \bK \rightarrow \bZ$ and $\up : \bK\rightarrow \bK'$.
We similarly define  the lifting map $\ell' : \bK' \rightarrow \bZ$ and $\down : \bK'\rightarrow \bK$ by
\begin{equation}\label{eq:down}\down (c)=\ell' (c)\bmod p      \text{ \,\,\, for \,\,} c \in \bK'.\end{equation}
These four maps can be extended to the corresponding group algebras by  coefficientwise application.
        Given $a$ and $b$ as in Equation
    (\ref{eq:ab}) we define
    \[A = \ell (a)= \sum_{\sigma \in G}\ell (a_\sigma)\sigma\text{ \,\,\, and \,\,   } B = \ell (b)=\sum_{\sigma \in G}\ell (b_\sigma)  \sigma  \text{\,\, in \,\, }  \bZ[G] \text{ \,\,\, and \,\,   } C=AB.\]
    The coefficients in $C$ belong to the interval $[0,\co . (p-1)^2 ]$. So \[C= \ell'( (A\mod p')\times (B\mod p')) \text{\,\,\, and \,\,\,} ab = \down (\up (a)\up (b)).\]
Using Proposition \ref{prop:enro} we deduce
\begin{proposition}\label{prop:primfi}
There exists an  absolute constant  $\cQ$ such that the following is true.
Let $G$, $\co$, $e$ be as
in Section \ref{sec:multiv}. Let  $\bK = \bZ/p\bZ$ be a prime field.
  There exists a prime integer $p'\leqslant \cQ (\co . p)^{11}$
  and a straight-line program of length smaller than
  $\cQ.\co.\log \co$ that computes
  the product   $c=\sum_g c_g[g]$ of two elements
  $a=\sum_g a_g[g]$ and $b=\sum_g b_g[g]$ in $\bK [G]$ given by their coefficients $(a_g)_g$ and
  $(b_g)_g$. The operations in this straight-line program are additions and multiplications in
  $\bZ/p'\bZ$ and evaluations
  of the maps $\up$ and $\down$ defined in Equations (\ref{eq:up}) and (\ref{eq:down}).
\end{proposition}

Now let $\bL$ be a field  extension of degree $d$ of $\bK=\bZ/p\bZ$.
We assume that elements in $\bL$ are represented by their
coordinates
in some $\bK$-basis of $\bL$.
The bilinear part of one  multiplication in $\bL [G]$
reduces to $\mu_p(d)$ multiplications in $\bK[G]$
where $\mu_p(d)$ is 
the $\bK$-bilinear complexity of multiplication in $\bL$.
Work by Chudnovsky \cite{CC}, Shparlinski, Tsfasmann, Vladut \cite{STV},
Shokrollahi \cite{SH}, Ballet  and Rolland \cite{BR,BAL},
Chaumine \cite{CHA},  Randriambololona \cite{RAN}
and others imply that  $\mu_p(d)$  is bounded by an absolute constant
times $d$.  
We deduce  the following theorem.

\begin{theorem}\label{prop:fifi}
  There exists an  absolute constant  $\cQ$ such that the following is true.
  Let $G$ be a  finite commutative group of order  $\co$ and exponent  $e$. 
  Let $\bK = \bZ/p\bZ$ and $\bL$  a field extension
  of degree $d$ of $\bK$.
  There exists a prime integer $p'\leqslant \cQ  (\co . p)^{11}$
  and a straight-line program of length $\leqslant \cQ (d.\co.\log \co+d^2.\co)$ that computes
  the product   $c=\sum_g c_g[g]$ of two elements
  $a=\sum_g a_g[g]$ and $b=\sum_g b_g[g]$ in $\bL [G]$ given by their coefficients $(a_g)_g$ and
  $(b_g)_g$. The operations in this straight-line program are additions and multiplications in   $\bZ/p\bZ$ and in 
  $\bZ/p'\bZ$ and evaluations
  of the maps $\up$ and $\down$ defined in Equations (\ref{eq:up}) and (\ref{eq:down}).
\end{theorem}

\begin{remark} The $d^2.\co$ summand in the complexity comes from the linear part in the Chudnovsky algorithm
for multiplication in  finite field extensions.\end{remark}

\section{Constructing functions in the Hilbert class field}\label{sec:constf}

We have defined in Section \ref{sec:comm} matrices $\cE$, $\cC$
and $\cI$ for the evaluation and interpolation of
global sections  of a $G$-equivariant
invertible sheaf on a curve $Y$ acted on freely by
a commutative group $G$. 
We have seen in Sections \ref{sec:comm}, \ref{sec:pade}, and \ref{sec:ga}
how to efficiently compute with these  matrices. In this section
we consider   the problem of computing these matrices.

We recall in Section \ref{sec:cftj} the necessary background from
class field theory of function fields over a finite field.
We illustrate the constructive  aspects of class fields
on  a small   example in section \ref{sec:excftj}.
An important feature of this method is that we only work
with divisors and functions on $X$, the quotient of $Y$ by $G$. This is of some
importance since in the applications presented in Sections \ref{sec:interpol} and \ref{sec:gc}  the genus
of $Y$ is much larger (e.g. exponentially) than the genus of $X$.

\subsection{Class field theory and the jacobian variety}\label{sec:cftj}

Let 
$X$ be a projective, smooth, absolutely integral  curve  over a finite field $\bK$ of characteristic $p$.
Let $\bKb$ be an algebraic closure of $\bK$.
We need  an abelian unramified cover $\tau : Y \rightarrow X$ over $\bK$, with
$Y$ 
absolutely integral.
We will require that $Y$ has a $\bK$-rational point $Q_1$. This
implies that $\tau$ is completely split above $P_1=\tau (Q_1)$.

According to class field theory \cite{gacc,rosen} there is a maximal abelian unramified
cover of $X$ over $\bK$ that splits totally above $P_1$. We briefly recall
its geometric construction.
Let $J_X$ be the jacobian variety of $X$ and let \[j_X : X\rightarrow J_X\]
be the Jacobi map with origin $P_1$.
Let \[F_\bK :  J_X\rightarrow J_X\] be the Frobenius endomorphism of degree
$|\bK|$, the cardinality of $\bK$.
The endomorphism  \[\wp = F_\bK-1 : J_X\rightarrow J_X\] is  an unramified Galois cover between $\bK$-varieties
with Galois group $J_X(\bK)$.
We denote by \[\taum : \Ym \rightarrow X\]
the pullback
of $\wp$ along $j_X$. This is the maximal
abelian unramified
cover of $X$ that splits totally above $P_1$.
Any such cover $\tau : Y \rightarrow X$
is thus a quotient of $\taum$ by some subgroup $H$
of $J_X(\bK)$. We set $G=J_X(\bK)/H$ and notice
that $G$ is at the same time
the fiber of $\tau$ above $P_1$ and its Galois
group, acting by translations in $J_X/H$.
\begin{equation*}
\begin{tikzcd}
  J_X(\bK)  \ar[r,hook] \ar[d]  & \Ym  \ar[r,hook] \ar[d] & J_X   \ar[d,"H"] \ar[dd,bend left = 60, "\wp"]  \\
  G=J_X(\bK)/H \ar[r,hook]  \ar[d] & Y  \ar[r,hook] \ar[d,"\tau "] &  J_X/H   \ar[d,"G"]\\
  0=P_1 \ar[r,hook] &X \ar[r,hook] &J_X &
\end{tikzcd}  
\end{equation*}
Let $P$ be a $\bK$-rational point on $X$ and let $\Qm$ be any point on $\Ym (\bKb)$ such
that
\[\taum (\Qm)=\wp(\Qm)=P.\]  We have $F_\bK(\Qm)=\Qm+P$.
So the Artin map and the Jacobi map coincide, and the decomposition group
of any place on $Y$ above $P$ is the subgroup of 
$G$ generated by $P$ itself.
In particular the fiber of $\tau$ above $P$ splits over $\bK$ if and only
if $P$ is sent into  $H$  by the Jacobi map. Equivalently the class of $P-P_1$
belongs to $H$.

\subsection{An example}\label{sec:excftj}

In this section $\bK$ is the field with three elements and $X$
is the plane projective curve with homogeneous equation
\[Y^2Z^3=X(X-Z)(X^3+X^2Z+2Z^3).\]
This is a smooth  absolutely integral curve of genus $2$.
The characteristic polynomial of the Frobenius of $X/\bK$
is \begin{equation}\label{eq:lpol}\chi_\bK(t)=t^4 + t^3 + 2t^2 + 3t + 9.\end{equation}
The characteristic polynomial of the Frobenius
of a curve over a finite field (given by a reasonable model) can be computed in time polynomial in $p.g.n$ where
$p$ is the characteristic of the field, $n$ its degree over the prime field, and
$g$ the genus of the curve, using  $p$-adic methods introduced by
Kato-Lubkin \cite{kl}, Satoh \cite{satoh}, Mestre
\cite{mestre}, Kedlaya \cite{ked}, Lauder and Wan \cite{lw}
and widely extended  since then.

When the  genus of the curve  is fixed, the characteristic polynomial of the Frobenius
can be computed in time polynomial in the logarithm of the cardinality
of $\bK$, using
 $\ell$-adic methods introduced by Schoof \cite{schoof} and generalized by Pila \cite{pila}.

We deduce from Equation (\ref{eq:lpol}) that the jacobian variety $J_X$
of $X$ has \[\chi_\bK (1)=16\] rational points. There are $5$ places of degree $1$ on $X$.
We let $P_1$ be  the unique place at $(0,1,0)$ and let \[P_2=(0,0,1), \,\, P_3=(1,0,1), \,\, P_4=(2,2,1), \,\, P_5=(2,1,1).\]

The Picard group $J_X(\bK)$ is the direct sum
of a  subgroup of order $8$
generated by the class of  $P_4-P_1$
and a subgroup of order $2$ generated by
$P_2-P_1$.
The class of $4(P_4-P_1)$
is the class of  $P_3-P_1$.
The classes of $P_2-P_1$
and $P_3-P_1$ generate a subgroup $H$ of $\Pic^0(X)$
isomorphic to $(\bZ/2\bZ)^2$. 
The quotient group \[G=J_X(\bK)/H=\Pic^0(X)/H\] is cyclic of order $4$ generated by $P_4-P_1$. 
So the subcover $\tau : Y \rightarrow X$ of $\Ym$ associated  with
$H$ is cyclic of order $4$. And the fibers above $P_1$, $P_2$, and $P_3$ in this cover
all split over $\bK$. We will work with this cover.

According to Kummer theory, there is a duality  (as group schemes)
between the prime to $p$ part of $\Pic^0(X)$
and the {\'e}tale part of the kernel of $F_\bK-p$. Associated to the quotient
$G = \Pic^0(X)/H$ there must be a  subgroup scheme isomorphic to $\mmu_4$
inside the latter kernel.

We let $\zeta$
be a primitive fourth root of unity
in $\bKb$ and 
denote by  $\bL$ the degree two extension of $\bK$ generated by $\zeta$.
In order to find the group of order $4$  we are interested in,
we use  algorithms to compute the kernels of $F_\bK-1$ and $F_\bK-p$
described in \cite[Chapter 13]{booktau}.
The idea is to pick   random elements in $J_X(\bL)$ and project
them onto the relevant characteristic subspaces for the action  of $F_\bK$, using
our  knowledge  of the characteristic polynomial $\chi_\bK$.
We set \[P_6=(2\zeta, 2) \text{\,\,\, and \,\,\,}  \Gamma = 2(P_6-P_4)\]
and find that the class $\gamma $ of $\Gamma$
is of order $4$ and satisfies 
\[F_\bK(\gamma)=3\gamma.\]
Thus $\gamma$
generates the group
  we were looking for.
There is a unique function $R$ in $\bL (X)$  with divisor $4\Gamma$ and taking value $1$ at $P_1$.
The cover $\tau : Y \rightarrow X$  we are interested in is obtained by adding a $4$-th root $r$ of
$R$ to $\bL (X)$. To be quite precise this construction produces the base change to $\bL$
of the cover we are interested in. This will be fine for our purpose. So we let
\[r=R^{1/4}\]
be the $4$-th root of $R$ taking value $1$ at $Q_1$. Equivalently we define $Q_1$ to be the point over $P_1$
where $r$ takes the value $1$.
With the notation of Section \ref{sec:orths} we take \[D=2P_5  \text{ \,\,\,  and  \,\,\, } P=P_1+P_2+P_3.\]
We let  $E$ be  the pullback of $D$
by $\tau$ and  $Q$ the pullback of $P$. 
We expect \[\cL (E)=H^0(Y, \cO_Y (E))\] to be a free  $\bK [G]$-module of rank \[\deg(D)-g_X+1=1.\]
This will be confirmed by our computations.
Because the fibers above $P_1$, $P_2$ and $P_3$ all split over $\bK$,
the evaluation map $\cL (E) \rightarrow \bA$ is described by a $3\times 1$
matrix
with coefficients in $\bK[G]$.

For every $2\leqslant i\leqslant 3$ we choose a $4$-th root of $R(P_i)$ in $\bL$. This amounts
to choosing a point $Q_{i,1}$ in the fiber of $\tau$ above $P_i$.
We let  $\sigma$  be the unique element in $G$ that sends $r$ to $\zeta.r$ so
\[G \owns \sigma : r \mapsto \zeta .r .\]

The $\bK$-vector space  $\cL (E)$ decomposes over $\bL$ as a sum of four eigenspaces  associated
to the four eigenvalues $1$, $\zeta $, $\zeta^2=-1$, $\zeta^3=-\zeta $ of $\sigma$.
Let $0\leqslant j\leqslant 3$ and let $f$ be an eigenfunction in  $\cL (E)$
associated with the eigenvalue $\zeta^j$. Then the quotient  $f/r^j$
is invariant by $G$ and its
divisor satisfies
\[(f/r^j) \geqslant -E - j.(r)= -E-j.\tau^*(\Gamma).\] So $f/r^j$ can be seen as a function
on $X$ with divisor bigger than or equal to  $-D -j\Gamma$.
The eigenspace $\cL (E)_j$ associated to $\zeta^j$ is thus obtained as the image of the map
\begin{equation*}
\xymatrix@R-2pc{
H^0(X, \cO_X(D+j\Gamma))  \ar@{->}[r]  & \cL (E)_j\\
 F  \ar@{|->}[r] & f= Fr^j}
\end{equation*}

Evaluating $f$ at $Q_{i,1}$ for $1\leqslant i\leqslant 3$ then reduces to evaluating $F=f/r^j$ at $P_i$
and multiplying the result by the chosen $4$-th root of $R(P_i)$,
raised to the power $j$.

This remark enables us to compute a $\bK$-basis of $\cL (E)$ consisting of  eigenfunctions of $\sigma$
and to evaluate the functions in this basis at the $(Q_{i,1})_{1\leqslant i\leqslant 3}$
without ever writing equations for $Y$. We only need to compute the Riemann-Roch spaces associated
to the divisors $D+j\Gamma$ on $X$ for $0\leqslant j\leqslant 3$. The Riemann-Roch space of  a divisor $D=D_+-D_{-}$
on a curve $X$ 
is computed in time polynomial in the genus of  $X$
and the degrees of the positive and negative parts
$D_+$ and $D_{-}$ of $D$,
using Brill-Noether algorithm and its many variants.
See \cite{HUIE, VOLC, HESS} and the most efficient general algorithm due
to Makdisi \cite{mak1, mak2}.
In case the exponent of $G$ is large, we may have to compute linear spaces
like $H^0(X,\cO_X(D+j\Gamma))$ for large $j$. In that case, one
should use the method introduced by Menezes, Okamoto, and Vanstone \cite{MOV}
in the context of pairing computation, in order  to replace $j$ by its
logarithm in the complexity. 

Passing from the values of the eigenfunctions to the evaluation matrix $\cE$ reduces to applying an inverse
Fourier transform. We find
\begin{equation*}\label{eq:cE}
  \cE = \begin{pmatrix} 1 \\ e_{1,2} \\ e_{1,3}\end{pmatrix} \text{\,\,\, with \,\,}  e_{1,1} = 1,
  \,\, e_{1,2}=1+2\sigma +2\sigma^2+2\sigma^3, \,\, e_{1,3} =
  2+2\sigma +2\sigma^2+\sigma^3.\end{equation*}
Having a unit for $e_{1,1}$ is quite convenient.
In general one says that
$\cE$ is systematic
when the top square submatrix
is the identity.
This is possible  when the first points $Q_{i,1}$
form a basis for the dual of  $\cL(E)$.
This situation is generic in some sense
 but not granted. From  a systematic matrix $\cE$  
it is trivial to deduce the associated checking and interpolation matrices

\begin{equation*}\label{eq:cC}
  \cC = \begin{pmatrix} e_{1,2} & e_{1,3}\\ -1 & 0 \\ 0&-1 \end{pmatrix} \text{\,\,\, and \,\,}  \cI = \begin{pmatrix}1 &0&0  \end{pmatrix}.\end{equation*}

\begin{remark}\label{rem}
  We may wonder how general is the method presented above.
  The approach via Kummer theory applies as long
  as the order $\co$ of $G$ is prime to $p$. In case the order
  of $G$ is a power of $p$, one may try to use Hasse-Witt theory
  instead, following the rather effective presentation in Serre \cite{Sertop}.
  When $\co$ is  neither prime to $p$ nor
  a power of $p$ we do not know any better method than the general
  purpose algorithm  in \cite{HESS}.
\end{remark}  
\section{Interpolation on algebraic curves}\label{sec:interpol}

In this section we recall two classical applications
of interpolation on algebraic curves over finite fields and
illustrate  the benefit of $\bK[G]$-module structures
in this context.
Section \ref{sec:multens} is concerned
with the multiplication tensor in finite fields.
In Sections \ref{sec:geocod} and \ref{sec:deco} we see that geometric
codes associated to $G$-equivariant divisors can be encoded
in quasi-linear time and decoded
in quasi-quadratic time if $G$ is
commutative, acts freely, and the code is long enough.

\subsection{The complexity of multiplication in finite fields}\label{sec:multens}

The idea of using Lagrange interpolation over an algebraic
curve to multiply two elements in a finite field 
is due to Chudnovsky \cite{CC} and has been developed by
Shparlinski, Tsfasmann and Vladut  \cite{STV}, Ballet
and Rolland \cite{BR}, Chaumine \cite{CHA},
Randriambololona \cite{RAN} and others.

Let $\bK$ be a finite field and let $\co \geqslant 2$ be an integer.
Let $Y$ be a  smooth, projective, absolutely
integral  curve over $\bK$ and $B$ a place of degree $\co$ on $Y$. Let  $\bL = H^0(B,  \cO_B)$ be the residue field
at $B$. We choose a divisor $E$ disjoint from $B$ and assume
that the evaluation map \[e_B : H^0(Y, \cO_Y(E))\rightarrow \bL\]
is surjective so that elements in $\bL$ can be represented
by functions in $H^0(Y, \cO_Y(E))$.
The latter functions will be characterized by their values at
a collection $(Q_i)_{1\leqslant i\leqslant N}$
of $\bK$-rational points on $Y$.
We denote by
\[e_Q : H^0(Y, \cO_Y(2E)) \rightarrow \bK^N\]
the evaluation map at these points which  we assume to be injective.
The multiplication of two elements $e_B(f_1)$
and $e_B(f_2)$ in $\bL$ can be achieved by evaluating $f_1$
and $f_2$ at the $Q_i$, then multiplying each
$f_1(Q_i)$ by the corresponding $f_2(Q_i)$, then
finding the unique function $f_3$
in $H^0(Y, \cO_Y(2E))$ taking value $f_1(Q_i)f_2(Q_i)$
at $Q_i$, then computing $e_B(f_3)$.
The number of bilinear multiplications in $\bK$
in the whole process is equal to $N$.

This method uses curves over $\bK$
with arbitrarily large genus 
having a number of $\bK$-points bigger than some positive constant
times their genus. It bounds the bilinear complexity of multiplication
in $\bL/\bK$ by an absolute constant times the degree $\co$
of $\bL$ over $\bK$,
but it says little abound the linear part of the algorithm, that is
 evaluation
of the maps $e_B$ and $e_Q$ and their right (resp. left) inverses.

Now assume that the group of $\bK$-automorphisms
of $Y$ contains  a cyclic subgroup
$G$ of order $\co$
acting freely on $Y$.
Let  $\tau : Y\rightarrow X$ be  the quotient
by $G$ map. Assume that $B$ is the fiber of $\tau$
above some rational point $a$ on $X$. Assume that $E$ (resp. $Q$) is the pullback
by $\tau$ of
a divisor $D$ (resp. $P$) on $X$.
Under mild conditions, all the linear spaces above become
free $\bK [G]$-modules and the evaluation maps are $G$-equivariant.
A computational consequence is that the linear part in the
Chudnovsky algorithm
becomes quasi-linear in the degree $\co$ of the extension $\bL/\bK$.
This remark has been exploited in \cite{CoLe, coez} to  
bound the complexity of multiplication of two elements
in a  finite field given by their coordinates
in a normal basis. The decompositions of the multiplication
tensor that are proven to exist in \cite{coez} can
be actually computed using the techniques
presented in  Section
\ref{sec:constf}.

\subsection{Geometric codes}\label{sec:geocod}

The construction of error correcting codes
by evaluating functions on 
algebraic curves of higher genus is due to Goppa \cite{gop1,go2}.
Let $Y$ be a smooth, projective, absolutely
integral curve over a finite field $\bK$ of characteristic $p$.
Let  $d$ be  the degree of $\bK$ over the prime field $\bZ/p\bZ$.
Let $g_Y$ be the genus of $Y$.
Let $Q_1$, \ldots, $Q_N$ be pairwise distinct
$\bK$-rational points on $Y$. Let $t_i$ be a uniformizing parameter at $Q_i$.
Let $E$ be a divisor that is disjoint  from $Q = Q_1+\dots +Q_N$.
Assume that \begin{equation}\label{eq:encad}
  2g_Y-1\leqslant \deg (E)\leqslant \deg (Q)-1.\end{equation}
Let \[\bA = H^0(Q, \cO_Q) =H^0(Y, \cO_Y/\cO_Y(-Q))\simeq \bK^N\] be the residue algebra at $Q$.
Let \[\hat\bA = H^0(Y, \Omega_{Y/\bK}(-Q)/\Omega_{Y/\bK}) \simeq  \bigoplus_{i=1}^N\bK \frac{dt_i}{t_i} \simeq \bK^N\]
be the dual of $\bA$.
Evaluation at the $Q_i$ defines an injective linear map \[\cL(E)=H^0(Y, \cO_Y(E))\rightarrow \bA.\]
We similarly define an injective linear map \[\itOmega(-Q+E)=H^0(Y, \Omega_{Y/\bK}(-Q+E))\rightarrow \hat\bA.\]
The two vector subspaces $\cL(E)$ and $\itOmega(-Q+E)$ are orthogonal to each other for the canonical duality pairing.
They can be considered as linear codes over $\bK$ and
denoted by $C_\cL$ and $C_\Omega$ respectively.
The code $C_\cL$ has
length $N$, dimension \[K=\deg (E)-g_Y+1\] and
minimum distance greater than or equal to
$N-\deg(E)$.
Given a basis of $\cL(E)$ one defines the generating  matrix $\cE_E$
of the code  $C_\cL$ to be
the $N\times K$-matrix of the injection $\cL(E)\rightarrow \bA=\bK^N$.
One similarly defines  the parity-check matrix $\cC_E$ to be the $N\times (N-K)$-matrix
of $\itOmega(-Q+E)\rightarrow \hat\bA$.
We finally  denote by  $\cI_E$ the $K\times N$-matrix of some projection of $\bA$ onto $C_\cL$.
A message of length $K$ is encoded  by multiplying the corresponding column on the left
by $\cE_E$. The received word is checked by multiplying it on the left by the transpose
of $\cC_E$. And the initial message is recovered from a correct codeword applying the interpolation matrix
$\cI_E$.
In full generality, coding,  testing and interpolating respectively require $2NK$, $2N(N-K)$ and $2KN$
operations in $\bK$.

Assume now that the
group of $\bK$-automorphisms of $Y$ contains
a finite  commutative subgroup
$G$ of order $\co$ acting freely on $Y$.
Let  $\tau : Y\rightarrow X$ be the quotient
by $G$ map. Assume that $\co$ divides $N$ and
let \[n=N/\co.\]
Assume that $Q$ is the pullback by $\tau$
of a divisor \[P=P_1+\dots+P_n\] on $X$.
Assume that $E$ is the pullback  of some divisor
$D$ on $X$. We are thus  in the situation of Section \ref{sec:comm}.
The code $C_\cL$ is a free $\bK[G]$-submodule of $\bA$
of rank \[k=K/\co\]
and
$\cC_\Omega$ is its orthogonal module
for the $\bK[G]$-bilinear form defined in Section \ref{sec:duaa}.

The matrices $\cE_E$, $\cC_E$, and $\cI_E$ can be seen as matrices with coefficients
in $\bK[G]$ of respective  sizes $n\times k$, $n\times (n-k)$, and $k\times n$.
Coding now requires $2nk$ operations in $\bK[G]$ rather
than $2NK$ operations in $\bK$.
According to Theorem \ref{prop:fifi}, each such operation requires less
than $\cQ .  d^2. \co . \log \co$ operations in $\bZ/p\bZ$ and
$\bZ/p'\bZ$ where $p'\leqslant \cQ . (\co . p)^{11}$ for some absolute constant $\cQ$.
The total cost of coding is thus
bounded by an absolute  constant times
\[\frac{NK}{\co^2}.d^2.\co.\log (\co) .  (\log p+\log \co)^2 = N.d^2.\log (\co) .  k.(\log p+\log \co)^2\]
elementary operations.

\begin{remark}\label{remcod} Assuming that $\log \co$ is bigger than $k$ times
  a positive constant, 
  the cost of coding
is quasi-linear in the length $N$ of the code.
The same holds  for parity-checking and interpolating.
Indeed the action of a large commutative group
$G$ provides a significant computational
advantage. We shall see in Section \ref{sec:gc}
that geometric class field theory produces examples
of free  commutative  group actions meeting this condition.\end{remark}

\subsection{Basic decoding}\label{sec:deco}

Assume that we are in the situation of the beginning of
Section \ref{sec:geocod},
and that  we have received a message $r$ in $\bA=\bK^N$.
Let $c$ be the closest codeword to $r$
in $C_\cL$ for the Hamming distance
in $\bK^N$. Write \[r=c+\epsilon \]
and  $\epsilon$ the error vector. Let $f$ be the unique function
in $\cL(E)$ such that $f=c\bmod Q$.
The support  of the error vector $\epsilon$
is the effective divisor $\Epsilon$ consisting of all
points $Q_i$ where $\epsilon$
is not-zero.
The degree of $\Epsilon$
is the
number of errors in $r$.

The principle  of the basic decoding  algorithm \cite{jlkh,svl}
is: if
$a_0$ is  a small degree function
vanishing at every point in the support  $\Epsilon$
then $a_0r=a_0c \mod Q$
is the residue modulo $Q$ of an algebraic
function $a_0f$ of not too large degree.
This function can  be recovered from its values at $Q$ if
$N$ is large enough.
More concretely we let $E_0$ be some auxiliary divisor on $Y$
with degree at least  $g_Y$ and set \[E_1=E+E_0.\] Let  $\cP$ be the subspace of $\cL (E_0)$
consisting of all $a_0$ such that  there exists $a_1$
in $\cL(E_1)$ with  $a_0r=a_1\bmod Q$. Non-zero
elements in $\cP$ are denominators for $r$ in the sense of Section
\ref{sec:pade}. 
We just saw that every function in $\cL (E_0)$
vanishing at every point in the support of $\epsilon$
belongs to $\cP$. 

Conversely if $a_0$ is in $\cP$ then $a_0r$ belongs
to  $\cL(E_1)$ modulo $Q$. But $a_0c$ belongs
to   $\cL(E_1)$ modulo $Q$ also because $a_0$ is in $\cL (E_0)$
modulo $Q$ and $c$ is in $\cL (E)$
modulo $Q$.
So $a_0(r-c)=a_0\epsilon$  belongs
to   $\cL(E_1)$ modulo $Q$.
There is a function in $\cL(E_1)$
that is $a_0\epsilon$ modulo
$Q$. This function has  $N-\deg (\Epsilon)$ zeros
and degree at most $\deg(E_1)=\deg (E)+\deg(E_0)$. If we assume
that
\begin{equation}\label{eq:basic}
  \deg (\Epsilon)\leqslant N-1 -\deg (E) -\deg (E_0)
\end{equation}
  then the latter function must be zero. So
  $a_0$ vanishes at $\Epsilon$.
Assuming Equation (\ref{eq:basic}) we thus have $\cP = \cL(E_0-\Epsilon)$.
Assuming further that \begin{equation}\label{eq:basic2}
   \deg (\Epsilon)\leqslant  \deg(E_0)-g\end{equation} this space
   is non-zero. Computing it is a matter of linear algebra
   and requires a constant times $N^3$ operations
   in $\bK$. 
   Given any non-zero element $a_0$ in $\cP$
   we denote by $A_0$ the divisor consisting of
   all $Q_i$ where $a_0$ vanishes. The degree of $A_0$
   is bounded by $\deg E_0$. The error $\epsilon$
   is an element  in $\bA$ with support contained in $A_0$
   and such that $r-\epsilon$ belongs to $C_\cL$.
   Finding $\epsilon$ is a linear problem in $\leqslant \deg E_0$
   unknows
   and $N-\deg (E)+g_Y-1$ equations. The solution
   is unique because the difference of two solutions
   is  in $C_\cL$ and has at least
   $N-\deg (E_0)$ zeros. And this is strictly greater than $\deg (E)$
   by Equation (\ref{eq:basic}).

   Combining Equations  (\ref{eq:basic}) and  (\ref{eq:basic2})
   we see that the basic decoding algorithm corrects up
   to $\dbasic$ errors where \begin{equation}\label{eq:dbasic}
     \dbasic = \frac{N-\deg (E)-1-g_Y}{2}.\end{equation}

Assume now that the
group of $\bK$-automorphisms of $Y$ contains
a finite  commutative subgroup
$G$ of order $\co$ acting freely on $Y$.
Let  $\tau : Y\rightarrow X$ be the quotient
by $G$ map. Assume that $\co$ divides $N$ and
let $n=N/\co$. Assume that $Q$ is the pullback by $\tau$
of a divisor \[P=P_1+\dots+P_n\] on $X$.
Assume that $E$ is the pullback  of some divisor
$D$ on $X$. Assume that $E_0$ is the pullback
of some divisor $D_0$ on $X$. Assume
that  $\cL(E_0)$ 
contains a free module of rank
$\deg (D_0)-g_X+1$ over $\bK[G]$. According to Proposition
\ref{prop:p'}, such an $E_0$ exists if 
the order $\co$ of $G$ is  prime to $p$.
According to Proposition
\ref{prop:pfin}, such an $E_0$ exists if
the order $\co$ of $G$ is  a power of  $p$, and
the cardinality $q$ of
$\bK$ is at least $4$, and the genus of $X$ is at least $2$.
Another sufficient condition if that $\deg (D_0)\geqslant 2g_X-1$.
According to Proposition \ref{prop:finddenom}
we can find a denominator $a_0$
at 
the expense of
$\cQ .  (\co . n . \log (\co . n))^2$ operations
in $\bK$ and
$\cQ .\co . n^3\log (\co . n)$ operations
in $\bK[G]$.
According to Theorem \ref{prop:fifi}, each
operation
in  $\bK[G]$
requires less than \[\cQ .  d^2.\co . \log (\co). (\log p+\log \co)^2\]
elementary operations.
The total cost of finding a denominator
is  thus
bounded by an absolute  constant times  \[N^2.n.d^2.\log^4 (\co.n.p)\]
elementary operation.

\begin{remark}\label{remdecod} Assuming that $\log \co$ is bigger than $n$ times
  a positive constant, 
  the cost of finding
a denominator 
is quasi-quadratic in the length $N$ of the code.
Once
 found a denominator, the error  can be  found at the same
 cost. \end{remark}


 \section{Good geometric
    codes with  quasi-linear
   encoding}\label{sec:gc}

 In this section we specialize the constructions presented
 in Sections \ref{sec:geocod}
 and \ref{sec:deco} using curves with many points and their
Hilbert  class fields.
 We quickly review in Section \ref{sec:artin}
 some standard  useful results and observations
 which we apply in Section \ref{sec:asym} to the construction
 of families of good geometric codes having  quasi-linear encoding
 and a quasi-quadratic decoder.
 Recall that a family of codes over a fixed
 alphabet is said to be good when
 the length tends to infinity while  both the rate and the
 minimum distance
 have a strictly positive $\liminf$.
 
 \subsection{Controlling the class group and the Artin map}\label{sec:artin}

 We keep the notation from Section \ref{sec:cftj}.
 In particular $P_1$ is a $\bK$-rational point on $X$ and
 \[j_X : X\rightarrow J_X\]
 is the Jacobi map with origin $P_1$.
 For the applications
 we have in mind we need some control on the $\bK$-rational points on $X$,
 on the
 group $\Pic^0 (X)$ and most importantly on the image of $X(\bK)$ in $\Pic^0 (X)$ by the Jacobi map.
  A typical advantageous situation would be:
 \begin{enumerate}\label{cond:3cond}
 \item $X$ has enough $\bK$-rational  points, that is a  fixed
   positive constant times its genus $g_X$,
 \item  a fixed positive proportion of these points are mapped by $j_X$ into a subgroup $H$,
 \item $H$ is not too large i.e. 
   the quotient $\log |H|/\log |\Pic^0(X)|$ is smaller than
   a fixed constant smaller than $1$.
 \end{enumerate}

 A range of geometric techniques relevant to that problem
 is presented in Serre's course \cite{serreratio}
 with the related motivation of constructing  curves
 with many points.
 One says that (a family of) curves over a fixed finite field
 of cardinality $q$
 have many points
 when the ratio of the number of rational points by the genus tends to
 $\sqrt q -1$. 
 Modular curves $X_0(N)$  have many points  over finite fields with $p^2$ elements, corresponding to supersingular moduli, as was noticed
 by Ihara \cite{ihara} and  by Tzfasman, Vladut, and Zink \cite{tvz}. These
 authors  also found families of Shimura curves
 having many points over fields with cardinality a square.
 Garcia and Stichtenoth
 \cite{garsti} constructed for every square $q$ an infinite   tower of algebraic curves over $\bF_q$ such that the quotient
 of the number of $\bF_q$-points  by the genus converges to $\sqrt q -1$, and the quotient of the genera of two consecutive
 curves converges to $q$.

 As for conditions (2) and (3) above,
 it is  noted in \cite[5.12.4]{serreratio} 
that  the images by $j_X$ of $P_2$, \dots, $P_{n}$ generate a subgroup
 $H$ with at most $n-1$ invariant factors. If the class group  $J_X(\bK)$ has  $I\geqslant n-1$ invariant factors then the size
 of the quotient $G$ is bigger than or equal to the product of the $I-(n-1)$ smallest invariant factors
 of   $J_X(\bK)$.

 Another favourable situation
 exploited in \cite{queb, nixi, vdg,guxi} is when   $\bK$ has a strict subfield
 $\bk$ and $X$ is defined over $\bk$  and $P_1$ is $\bk$-rational.
Then the Jacobi map sends the points
in $X(\bk)$ into the  subgroup $J_X(\bk)$ of $J_X(\bK)$. We will use
this remark in the next section.

 \subsection{A construction}\label{sec:asym}

Let $\bk$ be a finite field with characteristic $p$. Let $q$
be the cardinality of $\bk$. We assume that $q$ is a square.
We consider a family of curves $(X_k)_{k \geqslant 1}$ over $\bk$
having many points over $\bk$. For example
we may take $X_k$ to be the $k$-th curve in the Garcia-Stichtenoth tower associated with $q$.
We denote  by $g_{X}$  the genus of $X_k$. We omit the index  $k$
in the sequel because  there is no risk of confusion.
We denote by $n$ the number of $\bk$-rational points on $X$.
We denote these points by  $P_1$, \ldots, $P_n$
and let   $P$ be the effective divisor
sum of all these points.
We let $\bK$ be
a non-trivial extension of $\bk$.
We will assume that the  degree of $\bK$
over $\bk$ is $2$ because higher values seem to bring nothing 
but disadvantages.
We denote by $T$ the quotient
\[T=J_X(\bK)/J_X(\bk).\]
We denote by $T_p$ the $p$-Sylow subgroup of $T$. We denote
by $T_{p'}$ the complement subgroup of $T_p$ in $T$.
Let  $G$ be 
the bigger among $T_p$ and $T_{p'}$.
This is a quotient of $T$.
Let  $H$ be  the kernel of the composite map
$J_X(\bK)\rightarrow T=J_X(\bK)/J_X(\bk) \rightarrow G$.
Let $\co$ be the order of $G$.
We note that
\[\#  J_X(\bK)/J_X(\bk) \geqslant    \left( {q}-1 \right)^{2g_X}/ \left( \sqrt{q}+1 \right)^{2g_X}=
  \left( \sqrt{q}-1 \right)^{2g_X}\] so
\begin{equation}\label{eq:grex}
  \co\geqslant \sqrt{\sharp T}\geqslant \left( \sqrt{q}-1 \right)^{g_X}\end{equation}
grows exponentially in $g_X$ provided $q\geqslant 9$.
Also $G$ is  a $p$-group or a  $p'$-group.
We find ourselves in the situation of Section \ref{sec:cftj}.
Let   $\Ym$  be  the maximal unramified cover of $X$ over $\bK$
which is totally decomposed over $\bK$ above $P_1$. Let  $Y$
be the quotient of $\Ym$ by $H$. The fibers
of \[\tau : Y\rightarrow X\] above the points $P_1$, \ldots ,
$P_n$ all split over $\bK$.
Let  $Q$ be  the pullback of $P$ by $\tau$. This is a divisor
on $Y$ of  degree \[N=\co . n.\]
We choose a  real number $\ddaleth$ such that
\begin{equation}\label{eq:encad2}
  0 < \ddaleth < \frac{\sqrt q}{2} -2.\end{equation}  Our goal is to correct
up to $\ddaleth . \co. g_X$ errors.
Let $D$ be a divisor on $X$ that is disjoint from $P$ and such that 
\[\deg (D)= \lceil (\sqrt q - 2-2\ddaleth)g_X \rfloor \]
the closest integer to $(\sqrt q - 2-2\ddaleth)g_X$.
Let $E$ be the pullback of $D$ by $\tau$.
We deduce from Equation (\ref{eq:encad2}) that condition
(\ref{eq:encad}) is met at least asymptotically.
 From $X$, $Y$, $E$, and $Q$  the construction in Section \ref{sec:geocod}
 produces a code  $C_\cL$ over the field $\bK$
 with $q^2$ elements, having  length \[N= \co .  n \simeq (\sqrt q -1) . \co . g_X\]
 and  dimension
\[K = \co .  (\deg (D)-g_X+1)\simeq (\sqrt q -3- 2\ddaleth) . \co . g_X.\]
We set $k=K/\co$ and  deduce from Equation (\ref{eq:grex}) that
the $\liminf$ of
$(\log \co)/n$ and
$(\log \co)/k$ are strictily positive.
As explained in Remark \ref{remcod}, 
this implies that the code $C_\cL$ can be encoded
and parity-checked in quasi-linear deterministic
time in its length $N$,
and   decoded with the same complexity when there are no errors.
     Using the basic decoding algorithm as in Section \ref{sec:deco}
     one can decode in the presence of errors in quasi-quadratic
     probabilistic (Las Vegas) time
     up to the distance \[\dbasic =
     \frac{N-\deg (E)-1-g_Y}{2} \simeq  \ddaleth . \co . g_X\]
   defined by Equation (\ref{eq:dbasic}) as explained in Remark
   \ref{remdecod}. We denote by
   $\deltabasic$ the 
   relative distance $\dbasic/N$. The existence of a divisor $D_0$ with
   all the properties required in Section \ref{sec:deco} is granted because
   $G$ is either a $p$-group or a $p'$-group. So we can apply Proposition
   \ref{prop:p'} or Proposition \ref{prop:pfin} depending on the case.
     This finishes the proof of the  theorem below.

\begin{theorem}\label{th}
  Let $p$ be a prime integer
  and let $q$ be a power of $p$. Assume that $q$
  is a square and \begin{equation}\label{eq:25}q\geqslant 25.\end{equation} Let
  $\ddaleth$ be a real such that
  \begin{equation}\label{eq:daleth}
    0 < \ddaleth < \frac{\sqrt q}{2} -2.\end{equation}
  There exists
   a family of linear error correcting codes over
  the field with $q^2$ elements having length $N$  tending
  to infinity and such that
  \begin{enumerate}
      \item      the rate $R$
  satisfies
  \begin{equation*}\label{eq:rate}
    \lim R  = \frac{\sqrt q -3-2\ddaleth}{\sqrt q -1}\end{equation*}
  \item for each code there exists a straight-line program that encodes  in 
    quasi-linear time in the length $N$,
  \item for each code there exists a computation tree that decodes  in quasi-quadratic
probabilistic (Las Vegas)
  time
  in the length $N$ up to the  relative distance $\deltabasic$  and
  \[\lim \deltabasic = \frac{\ddaleth}{\sqrt q-1}.\]
  \end{enumerate}
  \end{theorem}       

\bigskip

\begin{remark}  The complexity statements in
  the theorem above are  non-uniform in the sense that they
  bound the complexity of coding and decoding assuming
  that the code is given  by its generating,  parity-check and interpolation  matrices
having  coefficients in the group algebra $\bK[G]$. The theorem claims nothing
 about the complexity of finding these matrices.
 The example detailed in section \ref{sec:excftj}
 and remark \ref{rem} suggest that this complexity could be quasi-quadratic
 in the length $N$  of the code. Proving such  a  complexity
 result would probably be quite heavy due to the relative sophistication
 of the methods used to find e.g.  the interesting torsion points  in the Picard group.\end{remark}


\begin{remark}  A  calculation similar to the one
in \cite[\S 7.3]{LACHAUD} shows that for any $q \geqslant 47^2$, 
some among the codes constructed above are  excellent 
in the sense that the accumulation
point $(2\deltabasic , R)$ stands above
the Varshamov-Gilbert limit for codes over the field
with $q^2$ elements. To our knowledge these are the first excellent codes that can be encoded in quasi-linear time and
decoded in quasi-quadratic time. Recall that Reed--Solomon  codes can be encoded and decoded in quasi-linear time
but cannot be said to be asymptotically good  because the length of the code is bounded by the size of the
alphabet.\end{remark}

\begin{remark}  We compare  fast basic decoding of the codes in   Theorem \ref{th} as explained in Section \ref{sec:deco} with  
the general purpose algorithm
of  Beelen, Rosenkilde, Solomatov \cite{beelen}. Using the latter, one can   decode up to half the Goppa designed minimum distance.
 Inequalities  (\ref{eq:25})  and (\ref{eq:daleth})
  are then replaced
  by
  \begin{equation*}\label{eq:16}q\geqslant 16
    \text{\,\,\,\,\,\, and \,\,\,\,\,\,}     0 < \ddaleth < \frac{\sqrt q -3}{2},\end{equation*}
  and the limit of the rate becomes 
\begin{equation*}
  \lim R  = \frac{\sqrt q -2-2\ddaleth}{\sqrt q -1}.\end{equation*}
However the complexity of decoding is
then of order
$\mu^{\omega-1}(N+g_Y)$ where $N$ is the length
of the code, $\mu$ is the gonality of
$Y$, and $\omega$ is the exponent in the complexity
of matrix multiplication. Curves with many points have large gonality.
In particular $\mu\geqslant N/(q^2+1)$ in our situation, so that
for fixed $q$, the complexity of this decoder 
is of order greater than $N^{\omega}$. It is known \cite{laser}
that $2\leqslant \omega <  2.37286$ but it  is not granted 
that $\omega =2$.
\end{remark}


\bibliographystyle{amsplain}

\bibliography{couveignes-gasnier}

\providecommand{\bysame}{\leavevmode\hbox to3em{\hrulefill}\thinspace}
\providecommand{\MR}{\relax\ifhmode\unskip\space\fi MR }
\providecommand{\MRhref}[2]{%
  \href{http://www.ams.org/mathscinet-getitem?mr=#1}{#2}
}
\providecommand{\href}[2]{#2}
\begin{thebibliography}{10}

\bibitem{laser}
Josh Alman and Virginia~Vassilevska Williams, \emph{A refined laser method and
  faster matrix multiplication}, Proceedings of the 2021 {ACM-SIAM} Symposium
  on Discrete Algorithms, {SODA} 2021, Virtual Conference, January 10 - 13,
  2021 (D{\'{a}}niel Marx, ed.), {SIAM}, 2021, pp.~522--539.

\bibitem{BLB}
S.~Ballet and D.~Le~Brigand, \emph{On the existence of non-special divisors of
  degree {$g$} and {$g-1$} in algebraic function fields over {$\mathbb F_q$}},
  J. Number Theory \textbf{116} (2006), no.~2, 293--310.

\bibitem{BR}
S.~Ballet and R.~Rolland, \emph{Multiplication algorithm in a finite field and
  tensor rank of the multiplication}, J. Algebra \textbf{272} (2004), no.~1,
  173--185.

\bibitem{BAL}
St\'{e}phane Ballet, \emph{Curves with many points and multiplication
  complexity in any extension of {${\bf F}_q$}}, Finite Fields Appl. \textbf{5}
  (1999), no.~4, 364--377.

\bibitem{beelen}
Peter Beelen, Johan Rosenkilde, and Grigory Solomatov, \emph{Fast decoding of
  {AG} codes}, 2022.

\bibitem{Blu}
Leo~I. Bluestein, \emph{A linear filtering approach to the computation of
  discrete {F}ourier transform}, IEEE Transactions on Audio and
  Electroacoustics \textbf{18} (1970), 451--455.

\bibitem{bost}
Alin Bostan, Fr{\'e}d{\'e}ric Chyzak, Marc Giusti, Romain Lebreton,
  Gr{\'e}goire Lecerf, Bruno Salvy, and \'Eric Schost, \emph{Algorithmes
  efficaces en calcul formel}, August 2017, 686 pages. \'Edition 1.0.

\bibitem{bouralg8}
N.~Bourbaki, \emph{\'{E}l\'{e}ments de math\'{e}matique. {A}lg\`ebre.
  {C}hapitre 8. {M}odules et anneaux semi-simples}, Springer, Berlin, 2012,
  Second revised edition of the 1958 edition.

\bibitem{CHA}
Jean Chaumine, \emph{Multiplication in small finite fields using elliptic
  curves}, Algebraic geometry and its applications, Ser. Number Theory Appl.,
  vol.~5, World Sci. Publ., Hackensack, NJ, 2008, pp.~343--350.

\bibitem{CC}
D.~V. Chudnovsky and G.~V. Chudnovsky, \emph{Algebraic complexities and
  algebraic curves over finite fields}, J. Complexity \textbf{4} (1988), no.~4,
  285--316.

\bibitem{coez}
Jean-Marc Couveignes and Tony Ezome, \emph{The equivariant complexity of
  multiplication in finite field extensions}, J. Algebra \textbf{622} (2023),
  694--720.

\bibitem{CoLe}
Jean-Marc Couveignes and Reynald Lercier, \emph{Elliptic periods for finite
  fields}, Finite Fields Appl. \textbf{15} (2009), no.~1, 1--22.

\bibitem{currei}
Charles~W. Curtis and Irving Reiner, \emph{Representation theory of finite
  groups and associative algebras}, Pure and Applied Mathematics, vol. Vol. XI,
  Interscience Publishers, New York-London, 1962.

\bibitem{booktau}
Bas Edixhoven and Jean-Marc Couveignes (eds.), \emph{Computational aspects of
  modular forms and {G}alois representations}, Annals of Mathematics Studies,
  vol. 176, Princeton University Press, Princeton, NJ, 2011.

\bibitem{garsti}
Arnaldo Garc\'{\i}a and Henning Stichtenoth, \emph{A tower of
  {A}rtin-{S}chreier extensions of function fields attaining the
  {D}rinfeld-{V}l{a}du{t} bound}, Invent. Math. \textbf{121} (1995), no.~1,
  211--222.

\bibitem{gop1}
V.~D. Goppa, \emph{Codes on algebraic curves}, Dokl. Akad. Nauk SSSR
  \textbf{259} (1981), no.~6, 1289--1290.

\bibitem{go2}
\bysame, \emph{Algebraic-geometric codes}, Izv. Akad. Nauk SSSR Ser. Mat.
  \textbf{46} (1982), no.~4, 762--781, 896.

\bibitem{gop2}
\bysame, \emph{Geometry and codes}, Mathematics and its Applications (Soviet
  Series), vol.~24, Kluwer Academic Publishers Group, Dordrecht, 1988.

\bibitem{guxi}
Venkatesan Guruswami and Chaoping Xing, \emph{Optimal rate list decoding over
  bounded alphabets using algebraic-geometric codes}, J. ACM \textbf{69}
  (2022), no.~2, Art. 10, 48.

\bibitem{hb}
D.~R. Heath-Brown, \emph{Zero-free regions for {D}irichlet {$L$}-functions, and
  the least prime in an arithmetic progression}, Proc. London Math. Soc. (3)
  \textbf{64} (1992), no.~2, 265--338.

\bibitem{HESS}
F.~Hess, \emph{Computing {R}iemann-{R}och spaces in algebraic function fields
  and related topics}, J. Symbolic Comput. \textbf{33} (2002), no.~4, 425--445.

\bibitem{HUIE}
Ming-Deh Huang and Doug Ierardi, \emph{Efficient algorithms for the
  {R}iemann-{R}och problem and for addition in the {J}acobian of a curve}, J.
  Symbolic Comput. \textbf{18} (1994), no.~6, 519--539.

\bibitem{ihara}
Yasutaka Ihara, \emph{Some remarks on the number of rational points of
  algebraic curves over finite fields}, J. Fac. Sci. Univ. Tokyo Sect. IA Math.
  \textbf{28} (1981), no.~3, 721--724 (1982).

\bibitem{jlkh}
J{\o}rn Justesen, Knud~J. Larsen, H.~Elbr{\o}nd Jensen, Allan Havemose, and Tom
  H{\o}holdt, \emph{Construction and decoding of a class of algebraic geometry
  codes}, IEEE Trans. Inform. Theory \textbf{35} (1989), no.~4, 811--821.

\bibitem{Kal}
Erich~L. Kaltofen and B.~David Saunders, \emph{On {W}iedemann's method of
  solving sparse linear systems}, Applied Algebra, Algebraic Algorithms and
  Error-Correcting Codes, 9th International Symposium, AAECC-9, New Orleans,
  LA, USA, October 7-11, 1991, Proceedings (Harold~F. Mattson, Teo Mora, and
  T.~R.~N. Rao, eds.), Lecture Notes in Computer Science, vol. 539, Springer,
  1991, pp.~29--38.

\bibitem{kl}
Goro~C. Kato and Saul Lubkin, \emph{Zeta matrices of elliptic curves}, J.
  Number Theory \textbf{15} (1982), no.~3, 318--330.

\bibitem{ked}
Kiran~S. Kedlaya, \emph{Counting points on hyperelliptic curves using
  {M}onsky-{W}ashnitzer cohomology}, J. Ramanujan Math. Soc. \textbf{16}
  (2001), no.~4, 323--338.

\bibitem{mak1}
Kamal Khuri-Makdisi, \emph{Linear algebra algorithms for divisors on an
  algebraic curve}, Math. Comp. \textbf{73} (2004), no.~245, 333--357.

\bibitem{mak2}
\bysame, \emph{Asymptotically fast group operations on {J}acobians of general
  curves}, Math. Comp. \textbf{76} (2007), no.~260, 2213--2239.

\bibitem{LACHAUD}
Gilles Lachaud, \emph{Les codes g\'{e}om\'{e}triques de {G}oppa},
  Ast\'{e}risque, vol. 133-134, 1986, Seminar Bourbaki, 1984/85, exp. 641,
  pp.~189--207.

\bibitem{lw}
Alan G.~B. Lauder and Daqing Wan, \emph{Counting points on varieties over
  finite fields of small characteristic}, Algorithmic number theory: lattices,
  number fields, curves and cryptography, Math. Sci. Res. Inst. Publ., vol.~44,
  Cambridge Univ. Press, Cambridge, 2008, pp.~579--612.

\bibitem{MOV}
Alfred~J. Menezes, Tatsuaki Okamoto, and Scott~A. Vanstone, \emph{Reducing
  elliptic curve logarithms to logarithms in a finite field}, IEEE Trans.
  Inform. Theory \textbf{39} (1993), no.~5, 1639--1646.

\bibitem{mestre}
J.-F. Mestre, \emph{Lettre adress{\'e}e {\`a} {G}audry et {H}arley},
  \url{https://webusers.imj-prg.fr/~jean-francois.mestre/}, december 2010.

\bibitem{nixi}
Harald Niederreiter and Chaoping Xing, \emph{A general method of constructing
  global function fields with many rational places}, Algorithmic number theory
  ({P}ortland, {OR}, 1998), Lecture Notes in Comput. Sci., vol. 1423, Springer,
  Berlin, 1998, pp.~555--566.

\bibitem{pila}
Jonathan~S. Pila, \emph{Frobenius maps of {A}belian varieties and finding roots
  of unity in finite fields}, ProQuest LLC, Ann Arbor, MI, 1988, Thesis
  (Ph.D.)--Stanford University.

\bibitem{queb}
Heinz-Georg Quebbemann, \emph{Cyclotomic {G}oppa codes}, IEEE Trans. Inform.
  Theory \textbf{34} (1988), no.~5, 1317--1320, Coding techniques and coding
  theory.

\bibitem{RSR}
Lawrence~R. Rabiner, Ronald~W. Schafer, and Charles~M. Rader, \emph{The chirp
  {$z$}-transform algorithm and its application}, Bell System Tech. J.
  \textbf{48} (1969), 1249--1292.

\bibitem{RAN}
Hugues Randriambololona, \emph{Bilinear complexity of algebras and the
  {C}hudnovsky-{C}hudnovsky interpolation method}, J. Complexity \textbf{28}
  (2012), no.~4, 489--517.

\bibitem{rosen}
Michael Rosen, \emph{The {H}ilbert class field in function fields}, Exposition.
  Math. \textbf{5} (1987), no.~4, 365--378.

\bibitem{satoh}
Takakazu Satoh, \emph{The canonical lift of an ordinary elliptic curve over a
  finite field and its point counting}, J. Ramanujan Math. Soc. \textbf{15}
  (2000), no.~4, 247--270.

\bibitem{schoof}
Ren\'{e} Schoof, \emph{Elliptic curves over finite fields and the computation
  of square roots mod {$p$}}, Math. Comp. \textbf{44} (1985), no.~170,
  483--494.

\bibitem{Sertop}
Jean-Pierre Serre, \emph{Sur la topologie des vari\'{e}t\'{e}s alg\'{e}briques
  en caract\'{e}ristique {$p$}}, Symposium internacional de topolog\'{\i}a
  algebraica {I}nternational symposium on algebraic topology, Universidad
  Nacional Aut\'{o}noma de M\'{e}xico and UNESCO, M\'{e}xico, 1958, pp.~24--53.

\bibitem{gacc}
\bysame, \emph{Algebraic groups and class fields}, Graduate Texts in
  Mathematics, vol. 117, Springer-Verlag, New York, 1988.

\bibitem{serreratio}
\bysame, \emph{Rational points on curves over finite fields}, Documents
  Math\'{e}matiques (Paris), vol.~18, Soci\'{e}t\'{e} Math\'{e}matique de
  France, Paris, 2020, With contributions by Everett Howe, Joseph Oesterl\'{e}
  and Christophe Ritzenthaler.

\bibitem{SH}
Mohammad~Amin Shokrollahi, \emph{Optimal algorithms for multiplication in
  certain finite fields using elliptic curves}, SIAM J. Comput. \textbf{21}
  (1992), no.~6, 1193--1198.

\bibitem{STV}
Igor~E. Shparlinski, Michael~A. Tsfasman, and Serge~G. Vladut, \emph{Curves
  with many points and multiplication in finite fields}, Coding theory and
  algebraic geometry ({L}uminy, 1991), Lecture Notes in Math., vol. 1518,
  Springer, Berlin, 1992, pp.~145--169.

\bibitem{svl}
Alexei~N. Skorobogatov and Sergei~G. Vl\u{a}du\c{t}, \emph{On the decoding of
  algebraic-geometric codes}, IEEE Trans. Inform. Theory \textbf{36} (1990),
  no.~5, 1051--1060.

\bibitem{sagemath}
{The Sage Developers}, \emph{{S}agemath, the {S}age {M}athematics {S}oftware
  {S}ystem ({V}ersion 9.4)}, 2021, {\tt https://www.sagemath.org}.

\bibitem{tvc}
M.~A. Tsfasman and S.~G. Vl\u{a}du\c{t}, \emph{Algebraic-geometric codes},
  Mathematics and its Applications (Soviet Series), vol.~58, Kluwer Academic
  Publishers Group, Dordrecht, 1991.

\bibitem{tvz}
M.~A. Tsfasman, S.~G. Vl\u{a}du\c{t}, and Th. Zink, \emph{Modular curves,
  {S}himura curves, and {G}oppa codes, better than {V}arshamov-{G}ilbert
  bound}, Math. Nachr. \textbf{109} (1982), 21--28.

\bibitem{vdg}
Gerard van~der Geer, \emph{Hunting for curves with many points}, Coding and
  cryptology, Lecture Notes in Comput. Sci., vol. 5557, Springer, Berlin, 2009,
  pp.~82--96.

\bibitem{VOLC}
Emil~J. Volcheck, \emph{Computing in the {J}acobian of a plane algebraic
  curve}, Algorithmic number theory ({I}thaca, {NY}, 1994), Lecture Notes in
  Comput. Sci., vol. 877, Springer, Berlin, 1994, pp.~221--233.

\bibitem{Wied}
Douglas~H. Wiedemann, \emph{Solving sparse linear equations over finite
  fields}, {IEEE} Trans. Inf. Theory \textbf{32} (1986), no.~1, 54--62.

\bibitem{xyl}
Triantafyllos Xylouris, \emph{On the least prime in an arithmetic progression
  and estimates for the zeros of {D}irichlet {$L$}-functions}, Acta Arith.
  \textbf{150} (2011), no.~1, 65--91.

\end{thebibliography}

\end{document}